\documentclass[draft]{amsart}

\usepackage{amsmath}

\usepackage{amssymb}

\usepackage{xypic}

\usepackage[all]{xy}

\usepackage{latexsym}

\usepackage{amsfonts}

\usepackage{graphicx}
\usepackage{epstopdf}
\usepackage{graphicx,epsfig,float}

\usepackage[latin1]{inputenc}

\usepackage[T1]{fontenc}

\numberwithin{equation}{section}

\newcommand{\A}{\mathcal{A}}

\newcommand{\K}{\mathcal{K}}

\newcommand{\T}{\mathcal{T}}

\newcommand{\U}{\mathcal{U}}

\renewcommand{\mod}{\mathrm{Mod}}

\newcommand{\psh}{\mathrm{Psh}}

\newcommand{\cov}{\mathrm{Cov}}

\newcommand{\coh}{\mathrm{Coh}}

\newcommand{\op}{\mathrm{Op}}

\newcommand{\rc}{\mathbb{R}\textrm{-}\mathrm{c}}

\newcommand{\R}{\mathbb{R}}

\newcommand{\N}{\mathbb{N}}

\renewcommand{\P}{\mathcal{P}}

\newcommand{\F}{\mathcal{F}}

\newcommand{\I}{\mathrm{I}}

\newcommand{\iso}{\stackrel{\sim}{\to}}

\newcommand{\ho}{\mathcal{H}\mathit{om}}

\newcommand{\Ho}{\mathrm{Hom}}

\newcommand{\id}{\mathrm{id}}

\newcommand{\coker}{\mathrm{coker}}

\renewcommand{\dim}{\textbf{Proof.}}

\newcommand{\pE}{\overset{\hspace{0.1cm}}{\dot{E}}}

\newcommand{\pU}{\overset{\hspace{0.1cm}}{\dot{U}}}

\newcommand{\pV}{\overset{\hspace{0.1cm}}{\dot{V}}}

\newcommand{\RP}{\mathbb{R}^{{\scriptscriptstyle{+}}}}

\newcommand{\imin}[1]{#1^{-1}}

\newcommand{\lind}[1]{\underset{#1}{\underrightarrow{\lim}}}

\newcommand{\Lind}{\underrightarrow{\lim}}  

\newcommand{\indl}[1]{\underset{#1}{``\underrightarrow{\mathrm{lim}}\mbox{''}}}

\newcommand{\lpro}[1]{\underset{#1}{\underleftarrow{\lim}}}

\newcommand{\exs}[3]{0 \to {#1} \to {#2} \to {#3} \to 0}

\newcommand{\lexs}[3]{0 \to {#1} \to {#2} \to {#3}}

\newcommand{\rexs}[3]{{#1} \to {#2} \to {#3} \to 0}

\newtheorem{teo}{Theorem}[subsection]

\newtheorem{df}[teo]{Definition}

\newtheorem{cor}[teo]{Corollary}

\newtheorem{oss}[teo]{Remark}

\newtheorem{prop}[teo]{Proposition}

\newtheorem{lem}[teo]{Lemma}

\newtheorem{es}[teo]{Example}

\newtheorem{nt}[teo]{Notation}

\usepackage{latexsym}
\begin{document}

\title{\bf{SHEAVES ON  $\T$-TOPOLOGIES}}

\date{}

\author{M\'{a}rio J. Edmundo}

\address{Universidade Aberta \& CMAF Universidade de Lisboa\\ Av. Prof. Gama Pinto 2\\ 1649-003 Lisboa, Portugal}

\email{edmundo@cii.fc.ul.pt}

\author{Luca Prelli }
\address{CMAF Universidade de Lisboa\\ Av. Prof. Gama Pinto 2\\ 1649-003 Lisboa, Portugal}

\thanks{The first author was supported by Funda\c{c}\~ao para a Ci\^encia e a Tecnologia, Financiamento Base 2008 - ISFL/1/209. The second author  is a member of the Gruppo Nazionale per l'Analisi Matematica, la Probabilit\`a e le loro Applicazioni (GNAMPA) of the Istituto Nazionale di Alta Matematica (INdAM) and was supported by Marie Curie grant PIEF-GA-2010-272021. This work is part of the FCT project PTDC/MAT/101740/2008. }


\maketitle

\thispagestyle{empty}
\begin{abstract}
The aim of this paper is to give a unifying description of various constructions of sites (subanalytic, semialgebraic, o-minimal) and consider the corresponding theory of sheaves.  The method used applies to a more general context and gives  new results in semialgebraic and o-minimal sheaf theory.
\end{abstract}

\tableofcontents

\section*{Introduction}

Sheaf theory in some tame contexts such as semi-algebraic geometry (\cite{De91}), subanalytic geometry (\cite{KS01,Pr1}) and o-minimal geometry (\cite{ejp}) has had recently different applications in various fields of mathematics such as model theory \cite{Be07,Be09,ep1}, analysis \cite{KS01,Mo1,Mo2,Pr07} and representation theory \cite{AG08,AG10,Pr}. Each one of the above theories is very useful for the mentioned applications but has some elements which are missing in the other ones: the aim of this paper is to give a unifying description of all these various constructions (subanalytic, semialgebraic, o-minimal) using a modification of the notion of $\T$-topology introduced by Kashiwara and Schapira in \cite{KS01}. \\

The idea is the following: on a topological space $X$ one chooses a subfamily $\T$ of  open subsets of $X$ satisfying some suitable hypothesis, and for each $U \in \T$ one defines the category of coverings of $U$ as the topological coverings $\{U_i\}_{i \in I} \subset \T$ of $U$ admitting a finite subcover. In this way one defines a site $X_\T$ and studies the category of sheaves on $X_\T$ (called $\mod(k_\T)$). This idea was already present in \cite{KS01}. However in \cite{KS01}, the space $X$ is assumed to be Hausdorff, locally compact and the elements of $\T$ are assumed to have finitely many connected components.\\

The exigence to treat in a unifying way all the previous constructions,  to   treat also some non Hausdorff cases (as conic subanalytic sheaves which are related to the extension of the Fourier-Sato transform \cite{Pr07}) and the non-standard setting which appears naturally in the o-minimal context (where the elements of $\T$  are totally disconnected and never locally compact), motivates a modification of the definition of \cite{KS01}. In particular, in our definition we replace ``connectedness'' by the notion of $\T$-connectedness (which in the standard o-minimal context is connectedness). Remark that there are many important o-minimal expansions
$${\mathcal M}=({\mathbb R},<, 0,1,+,\cdot, (f)_{f\in {\mathcal F}} )$$
of the ordered  field of real numbers. For example ${\mathbb R}_{{\rm an} }$, ${\mathbb R}_{{\rm exp }}, $ ${\mathbb R}_{{\rm an},\, {\rm exp}}, $ ${\mathbb R}_{{\rm an} ^*}, $ $ {\mathbb R}_{{\rm an} ^*,\,{\rm exp}}$ see resp., \cite{dd,w,dm94,ds1,ds2}. For each such we have $2^{\kappa }$ many non-isomorphic non standard o-minimal models for each $\kappa >$ cardinality of the language. There is however a non-standard o-minimal structure
$${\mathcal M}=\left(\bigcup _{n\in {\mathbb N}} {\mathbb R}((t^{\frac{1}{n}})), <,  0,1, +, \cdot, (f_p)_{p\in {\mathbb R}[[\zeta _1, \ldots ,\zeta _n]]}\right)$$
which does not came from a standard one (\cite{LpRob06,HrPet07}).\\


With this more general notion of $\T$-space $X$  we study  the category of sheaves on the site $X_\T$. 
The natural functor of sites $\rho:X \to X_\T$ induces relations between the categories of sheaves on $X$ and $X_\T$, given by the functors $\rho_*$ and $\imin \rho$. The functor $\rho _*$ is fully faithful. Moreover when $X$ is locally weakly quasi-compact  there is a right adjoint to the functor $\imin \rho$, denoted by $\rho_!$. The functor $\rho_!$ is exact, commutes with $\Lind$ and $\otimes$ and is fully faithful. We introduce the category of $\T$-flabby sheaves (known as $sa$-flabby in \cite{De91} and as quasi-injective in \cite{Pr1}): $F \in \mod(k_\T)$ is $\T$-flabby if the restriction $\Gamma(U;F) \to \Gamma(V;F)$ is surjective for each $U,V \in \T$ with $U \supseteq V$. We prove that $\T$-flabby sheaves are stable under $\Lind$ and $\otimes$ and are acyclic with respect to the functor $\Gamma(U;\bullet)$, for $U \in \T$. More generally,
if one introduces the category $\coh(\T) \subset \mod(k_X)$ of coherent sheaves (i.e. sheaves admitting a finite resolution consisting of finite sums of $k_{U_i}$, $U_i \in \T$), then
$\T$-flabby sheaves are acyclic with respect to $\Ho_{k_\T}(\rho_*G,\bullet)$, for $G \in \coh(\T)$. Coherent sheaves also give a description of sheaves on $X_\T$: for each $F \in \mod(k_\T)$ there exists a filtrant inductive family $\{F_i\}_{i\in I}$ such that $F \simeq \lind i \rho_*F_i$. In fact, we have an equivalence between the categories $\mod(k_\T)$ and ${\rm Ind}(\coh(\T))$ the indization of the category $\coh(\T)$.\\

All of the above results and methods are new in the o-minimal context and most of them are new even in the semialgebraic case as well. On the other hand, we also introduce a method for studying  the category $\mod(k_\T)$ of sheaves on $\T$-spaces which is the fundamental tool in the semialgebraic and o-minimal case, namely,  we prove that as in \cite{ejp} the category of sheaves on $X_\T$ is equivalent to the category of sheaves on a locally  quasi-compact  space $\widetilde{X}_\T$, the $\T$-spectrum of $X$, which generalizes the notion of o-minimal spectrum as well as the real spectrum of commutative rings from real algebraic geometry. In particular, sheaves on the subanalytic site are sheaves on the $\T$-spectrum associated to the family of relatively compact subanalytic subsets. Such a result was not present in \cite{KS01}.\\ 

This theory can then be specialized to each of the examples we mentioned above: when $\T$ is the category of semialgebraic open subsets of a locally semialgebraic space $X$ we obtain the constructions (and the generalizations) of results of \cite{De91}, in particular, when $X$ is a Nash manifold, we recover the setting of \cite{Pr}. When $\T$ is the category of relatively compact subanalytic open subsets of a real analytic manifold $X$ we obtain the constructions and results of \cite{KS01,Pr1}.  Moreover, when $\T$ is the category of conic subanalytic open subsets of a real analytic manifold $X$ we obtain a suitable category of conic subanalytic sheaves considered in \cite{Pr07}. Finally, when $\T$ is the category of definable open subsets of a locally definable space $X$ we obtain in the definable case  the constructions of \cite{ejp} and we obtain new results in the o-minimal context generalizing those of the two previous cases.\\

The paper is organized in the following way. In Section \ref{Sec sheaves on lwqcs} we introduce the locally weakly quasi-compact spaces and study some properties of sheaves on such spaces. The results  of this section will be used in two crucial ways on the theory of sheaves on $\T$-spaces, they are required to show that: (i) when a $\T$-space $X$ is locally weakly quasi-compact, then   there is a right adjoint $\rho_!$ to the functor $\imin \rho$ induced by  the natural functor of sites $\rho:X \to X_\T$; (ii)  for a $\T$-space $X$, the category of sheaves on $X_\T$ is equivalent to the category of sheaves on a locally  quasi-compact  space $\widetilde{X}_\T$, the $\T$-spectrum of $X$. In Section \ref{Sec tspaces} we introduce the $\T$-spaces and develop the theory of sheaves on such spaces as already described above. 

\section{Sheaves on locally weakly quasi-compact spaces}\label{Sec sheaves on lwqcs}

Let $X$ be a non necessarily Hausdorff topological space. One denotes by $\op(X)$ the category whose
objects are the open subsets of $X$ and the morphisms are the inclusions. In this section we generalize some classical results about sheaves on locally compact spaces. For classical sheaf theory our basic reference is \cite{KS90}. We refer to \cite{Ta94} for an introduction to sheaves on Grothendieck topologies.

\subsection{Locally weakly quasi-compact spaces}

\begin{df} An open subset $U$ of $X$ is said to be relatively weakly quasi-compact in $X$ if, for any covering $\{U_i\}_{i \in I}$ of $X$, there
  exists $J \subset I$ finite, such that
$U \subset \bigcup_{i \in J} U_i.$\\
\end{df}


We will write for short $U \subset\subset X$ to say that $U$ is a relatively
weakly quasi-compact open set in $X$, and we will call $\op^c(U)$ the subcategory
of $\op(U)$ consisting of open sets $V \subset\subset U$. Note
that, given $V,W \in \op^c(U)$, then $V \cup W \in \op^c(U)$.



\begin{df}
A topological space $X$ is locally weakly quasi-compact if satisfies the following hypothesis for every $U,V\in \op(X)$
\begin{itemize}
\item[LWC1.]
Every  $x \in U$ has a fundamental neighborhood system $\{V_i\}$ with $V_i \in \op^c(U)$.
\item[LWC2.]
For every $U' \in \op^c(U)$ and $V' \in \op^c(V)$ one has $U'  \cap V' \in \op^c(U \cap V)$.
\item[LWC3.]
For every $U' \in \op^c(U)$ there exists $W \in \op^c(U)$ such that $U' \subset\subset W$.
\end{itemize}
\end{df}
Of course an open subset $U$ of a locally weakly quasi-compact space $X$ is also a locally weakly quasi-compact space. Let us consider some examples of locally weakly quasi-compact
spaces:

\begin{es}
{\em
A locally compact topological space $X$ is a locally weakly quasi-compact.
In this case, for $U,V \in \op(X)$ we have $V\subset \subset U$ if and only if $V$ is relatively compact subset of  $U$.
}
\end{es}

\begin{es}
{\em
Let  $X$ be a topological space with a basis of quasi-compact (i.e. each open covering admits a finite subcover) open subsets  closed under taking finite intersections. Then  $X$ is locally weakly quasi-compact and, for $U,V \in \op(X)$ we have $V\subset \subset U$ if and only if $V$ is contained in a quasi-compact subset of $U$. In this situation we have the following particular cases:
\begin{itemize}
\item[(i)]
$X$ is a Noetherian topological space (each open subset of $X$ is quasi-compact). This includes in particular: (a)  algebraic varieties over algebraically closed fields; (b) complex varieties (reduced, irreducible complex analytic spaces) with the Zariski topology.
\item[(ii)]
$X$ is a spectral topological space (in addition: (i) $X$ is quasi-compact;  (ii) $T_0$;   (iii) every irreducible closed subset is the closure of a unique point). This includes in particular: (a) real algebraic varieties over real closed fields; (b) the o-minimal spectrum of a definable space in some o-minimal structure.
\item[(iii)]
$X$ is an increasing union of open spectral topological spaces $X_i$'s, i.e. $X$ is the space $\bigcup _{i\in I}X_i$. This space $X$ has  a basis of quasi-compact open subsets  closed under taking finite intersections and in addition is: (i) not quasi-compact in general unless  $I$ is finite; (ii) $T_0$. This includes in particular: (a) the semialgebraic spectrum of locally semialgebraic space; (b) more generally, the o-minimal spectrum of a locally definable space in some o-minimal structure.
\end{itemize}
}
\end{es}

\begin{es}
{\em
Let $E$ be a real vector bundle over a locally compact space $Z$ endowed with the
 natural action $\mu$ of $\RP$ (the multiplication on the fibers). Let $\pE=E\setminus Z$,
  and for $U \in \op(E)$ set $U_Z=U\cap Z$ and $\pU=U \cap \pE$. Let
$E_{\RP}$ denote the space $E$ endowed with the conic topology
i.e. open sets of $E_{\RP}$ are open sets of $E$ which are
$\mu$-invariant. With this topology $E_{\RP}$ is a locally
weakly quasi-compact space and,  for $U,V \in \op(E_{\RP})$ we have  $V
\subset\subset U$ if and only if $V_Z \subset \subset U_Z$ in $Z$
and $\pV \subset\subset \pU$ in $\pE_{\RP}$ (the later is $\pE$ with the induced conic topology).
}
\end{es}

\subsection{Sheaves on locally weakly quasi-compact spaces}

Recall that $X$ is a non necessarily Hausdorff topological space.

\begin{df} Let $\U=\{U_i\}_{i \in I}$ and $\U'=\{U'_j\}_{j \in J}$
be two families of open subsets of $X$. One says that $\U'$ is a refinement of $\U$ if for each $U_i \in \U$ there is $U'_j \in \U'$ with $U'_j \subseteq U_i$.
\end{df}

One denotes by $\cov(U)$ the category whose objects are the
coverings of $U \in \op(X)$ and the morphisms are the refinements,
and by $\cov^f(U)$ its full subcategory consisting of finite
coverings of $U$.\\

Given $V \in \op(U)$ and $S \in \cov(U)$, one sets $S \cap V=\{U
\cap V\}_{U \in S} \in \cov(V)$.

\begin{df} The site $X^f$ on the topological space $X$ is the category $\op(X)$ endowed with
the fol\-lowing topology: $S \subset \op(U)$ is a covering of $U$
if and only if it has a refinement $S^f \in \cov^f(U)$.
\end{df}


\begin{df} Let $U,V \in \op(X)$ with $V \subset U$. Given $S=\{U_i\}_{i
\in I} \in \cov(U)$ and $T=\{V_j\}_{j\in J} \in \cov(V)$, we write
$T \subset\subset S$ if $T$ is a refinement of $S \cap V$, and
$V_j \subset U_i$ if and only if $V_j \subset\subset U_i$.
\end{df}



Let us recall the definitions of presheaf and sheaf on a site.

\begin{df} A presheaf of $k$-modules on $X$ is a functor
$\op(X)^{op} \to \mod(k)$. A morphism of presheaves is a morphism
of such functors. One denotes by $\psh(k_X)$ the category of
presheaves of $k$-modules on $X$.
\end{df}

Let $F \in \psh(k_X)$, and let $S \in \cov(U)$. One sets
$$F(S)=\ker \Big(\prod_{W \in S}F(W)
\rightrightarrows \prod_{W',W'' \in S}F(W' \cap W'')\Big).$$

\begin{df} A presheaf $F$ is separated (resp. is a
sheaf) if for any $U \in \op(X)$ and for any $S \in \cov(U)$ the
natural morphism $F(U) \to F(S)$ is  a monomorphism (resp. an
isomorphism). One denotes by $\mod(k_X)$ the category of sheaves
of $k$-modules on $X$.
\end{df}

Let $F \in \psh(k_X)$, one defines the presheaf $F^+$ by setting
$$F^+(U) = \lind {S \in \cov(U)} F(S).$$
One can show that $F^+$ is a separated presheaf and if $F$ is a
separated presheaf, then $F^+$ is a sheaf. Let $F \in \psh(k_X)$,
the sheaf $F^{++}$ is called the sheaf associated to the presheaf
$F$.\\


\begin{lem}\label{fac+} For $F \in \psh(k_X)$, and let $U \in \op(X)$. If $F$ is a sheaf on $X^f$, then  for any $V \in
\op^c(U)$ the morphism
\begin{equation}\label{factors+}
F^+(U) \to F^+(V)
\end{equation}
 factors through $F(V)$.
\end{lem}
\dim\ \ Let $S \in \cov(U)$, and set $S \cap V = \{W \cap V\}_{W
\in S}$. Since $V \in \op^c(U)$, there is a finite refinement $T^f \in
\cov^f(V)$ of $S \cap V$. Then the morphism (\ref{factors+}) is
defined by
\begin{eqnarray*}
F^+(U) & \simeq & \lind {S \in \cov(U)} \hspace{-3mm} F(S)\\
 & \to & \lind {S \in \cov(U)} \hspace{-3mm}F(S \cap V )\\
& \to & \lind {T^f \in \cov^f(V)}\hspace{-3mm} F(T^f)\\
& \to & \lind {T \in \cov(V)}\hspace{-3mm}F(T)\\
& \simeq & F^+(V).
\end{eqnarray*}
The result follows because $F(T^f) \simeq F(V)$. \qed

\begin{cor}\label{fac|} With the hypothesis of Lemma \ref{fac+}, we consider two co\-ve\-rings $S \in \cov(U)$ and $T \in
\cov(V)$. If $T \subset\subset S$, then the morphism
\begin{equation}\label{factorscov+}
F^+(S) \to F^+(T) \end{equation}
 factors through $F(T)$. In particular, if $T$ is finite, then the
 morphism \eqref{factorscov+} factors through $F(V)$.
\end{cor}

From now on we will assume the following hypothesis:

 \begin{equation}\label{hyplqc}
 \text{the topological space $X$ is locally weakly quasi-compact.}
 \end{equation}


\begin{lem}\label{TccS} Let $U \in \op(X)$, and consider a subset
$V \subset\subset U$. Then for any $S^f \in \cov^f(U)$ there
exists $T^f \in \cov^f(V)$ with $T^f \subset\subset S^f$.
\end{lem}
\dim\ \ Let $S^f=\{U_i\}$. For each $x \in U$ and $U_i \ni x$,
consider a $V_{x,i} \in \op^c(U_i)$ containing $x$. Set
$V_x=\bigcap_i V_{x,i}$, the family $\{V_x\}$ forms a covering of
$U$. Then there exists a finite subfamily $\{V_j\}$ containing
$V$. By construction $V_j \cap V \subset\subset U_i$ whenever $V_j
\subset U_i$. \qed

\begin{lem}\label{fac++} Let  $F \in \psh(k_X)$, and let $U \in \op(X)$. If $F$ is a sheaf on $X^f$, then for any $V \in
\op^c(U)$ the morphism
\begin{equation}\label{factors++}
F^{++}(U) \to F^{++}(V)
\end{equation}
 factors through $F(V)$.
\end{lem}
\dim\ \ Since $X$ is locally weakly quasi-compact, there exists $W \in
\op^c(U)$ with $V \subset\subset W$. As in Lemma \ref{fac+} we
obtain a diagram
$$
\xymatrix{F^{++}(U) \ar[r] \ar@ <-2pt> [d] & F^{++}(W) \ar[r]
\ar[d]
& F^{++}(V) \\
 \hspace{-15mm}\lind {S^f \in \cov^f(W)}\hspace{-3mm}
F^+(S^f) \ar@ <7pt> [r] \ar[ur] & \hspace{-4mm}\lind {T^f \in
\cov^f(V)}\hspace{-3mm} F^+(T^f). \ar[ur]. & }
$$

Since $X$ is locally weakly quasi-compact then by Lemma \ref{TccS} for any $S^f \in
\cov^f(W)$ there exists $T^f \in \cov^f(V)$ with $T^f
\subset\subset S^f$. By Corollary \ref{fac|} the morphism
$$F^+(S^f) \to F^+(T^f)$$
 factors through $F(T^f) \simeq F(V)$. Then the morphism
$$\lind {S^f \in \cov^f(W)}\hspace{-3mm}F^+(S^f) \to \hspace{-3mm}\lind {T^f \in
\cov^f(V)}\hspace{-3mm}F^+(T^f)$$ factors through $F(V)$ and the
result follows. \qed

\begin{cor}\label{faclimlim} Let  $F \in \psh(k_X)$. If $F$ is a sheaf on $X^f$, then:
\begin{itemize}
\item[(i)] for any $V \in \op^c(X)$ one has the isomorphism $\lind {U \supset\supset V}F(U)
\iso \lind {U \supset\supset V}F^{++}(U)$. \\
\item[(ii)] for any $U \in \op(X)$ one has the isomorphism $\lpro {V \subset\subset U}F(V)
\iso \lpro {V \subset\subset U}F^{++}(V)$.
\end{itemize}
\end{cor}
\dim\ \ (i) By Lemma \ref{fac++} for each $U \in \op(X)$ with $U
\supset\supset V$ we have a commutative diagram
$$
\xymatrix{ F^{++}(U) \ar[r]  \ar[rd] & F^{++}(V)  \\
F(U) \ar[u] \ar[r] & F(V) \ar[u] .}
$$
This implies that the identity morphism of $\lind {U
\supset\supset V}F(U)$ factors through $\lind {U \supset\supset V}
F^{++}(U)$. On the other hand this also implies that the identity
morphism of $\lind {U \supset\supset V}F^{++}(U)$ factors through
$\lind {U \supset\supset V} F(U)$. Then $\lind {U \supset\supset
V}F(U) \iso \lind {U \supset\supset V}F^{++}(U)$.

The proof of (ii) is similar. \qed

\begin{cor}\label{facKK} Let $X$ be a quasi-compact and locally weakly quasi-compact
space, and let $F \in \psh(k_X)$. If $F$ is a sheaf on $X^f$, then
 the natural morphism
\begin{equation}\label{factorsK}
F(X) \to F^{++}(X)
\end{equation}
is an isomorphism.
\end{cor}
\dim\ \ It follows immediately from Corollary \ref{faclimlim} (i)
with $V=X$. \qed \\

Let $\{F_i\}_{i \in I}$ be a filtrant inductive system in
$\mod(k_X)$. One sets
\begin{equation*}
\begin{array}{l}
\text{$\indl i F_i$ = inductive limit in the category of presheaves,} \\
\text{$\lind i F_i$ = inductive limit in the category of sheaves.}
\end{array}
\end{equation*}
Recall that $\lind i F_i=(\indl i F_i)^{++}$.

\begin{prop}\label{fac} Let $\{F_i\}_{i \in \I}$ be a filtrant inductive system in $\mod(k_X)$ and let $U \in \op(X)$.
Then for any $V \in \op^c(U)$ the morphism
\begin{equation*}
\Gamma(U;\lind i F_i) \to \Gamma(V;\lind i F_i)
\end{equation*}
 factors through $\lind i \Gamma(V;F_i)$.
\end{prop}
\dim\ \ By Lemma \ref{fac++} it is enough to show that $\indl i
F_i$ is a sheaf on $X^f$. Let $U \in \op(X)$ and $S \in
\cov^f(U)$. Since $\lind i$ commutes with finite projective limits
we obtain the isomorphism $(\indl i F_i)(S) \simeq \lind i
F_i(S)$. The result follows because $F_i \in \mod(k_X)$ for each
$i \in I$. \qed

\begin{cor}\label{faclim} Let $\{F_i\}_{i \in \I}$ be a filtrant inductive system in $\mod(k_X)$.
\begin{itemize}
\item[(i)] For any $V \in \op^c(X)$ one has the isomorphism $\lind {U
\supset\supset V,i}\Gamma(U;F_i) \iso \lind {U \supset\supset
V}\Gamma( U; \lind i F_i)$. \\
\item[(ii)] For any $U \in \op(X)$ one has the isomorphism $\lpro {V \subset\subset
U}\lind i \Gamma(V;F_i) \iso \lpro {V \subset\subset
U}\Gamma(V;\lind i F_i)$.
\end{itemize}
\end{cor}
\dim\ \ It follows from Corollary \ref{faclimlim} with $F=\indl i
F_i$. \qed

\begin{cor}\label{facK} Let $X$ be a quasi-compact and locally weakly quasi-compact
space. Then the natural morphism
$$ \lind i \Gamma(X;F_i) \to \Gamma(X;\lind i F_i)$$
is an isomorphism.
\end{cor}
\dim\ \ It follows from Corollary \ref{facKK} with $F=\indl i
F_i$. \qed

\begin{es}\label{exiUVi} Let us consider the formula
\begin{equation}\label{iUVi}
\lind {U \supset\supset V,i}\Gamma(U;F_i) \iso \lind {U
\supset\supset V}\Gamma( U; \lind i F_i)
\end{equation}
\begin{itemize}
\item[(i)] Let $X$ be a Noetherian space and let $V \in \op(X)$. Then $\Gamma(V;F) \simeq \lind {U \supset\supset
V}\Gamma(U;F)$, since every open set is quasi-compact and
\eqref{iUVi} becomes $\lind i \Gamma(V;F_i) \simeq \Gamma(V;\lind
i F_i)$.
\item[(ii)] Assume that $X$ has a basis of quasi-compact open subsets and let $V \in \op^c(X)$. Then $V$ is contained in a quasi-compact open subset of $X$ and $\lind {U \supset\supset V}\Gamma(U;F) \simeq  \lind {W \supset V}\Gamma(W;F)$, where $W$ ranges through the family of quasi-compact subsets of $X$.
\item[(iii)] Let $X$ be a locally compact space and let $V \in \op^c(X)$. Then $\Gamma(\overline{V};F) \simeq \lind {U \supset\supset
V}\Gamma(U;F)$, and \eqref{iUVi} becomes $\lind
i\Gamma(\overline{V};F_i) \simeq \Gamma(\overline{V};\lind
iF_i)$.
\item[(iv)] Let $E_{\RP}$ be a vector bundle endowed with the
conic topology, and let $V \in \op^c(E_{\RP})$. Then $\lind {U
\supset\supset V}\Gamma(U;F) \simeq \Gamma(K;F)$, where $K$ is the
union of the closures of $V_Z$ in $Z$ and $\pV$ in $\pE_{\RP}$,
and \eqref{iUVi} becomes $\lind i\Gamma(K;F_i) \simeq
\Gamma(K;\lind iF_i)$.
\end{itemize}
\end{es}

\begin{lem}\label{lproc} Let $F \in \psh(k_X)$. Then we have the isomorphism
$$
\lpro {V\subset\subset X} \lind {V\subset\subset W} F(W) \iso \lpro {V \subset\subset X} F(V).
$$
\end{lem}
\dim\ \ The result follows since for each $V \in \op^c(X)$ there exists $W \in \op^c(X)$ such that $V \subset\subset W$ since $X$ is locally weakly compact. Let $U,V \subset\subset X$ such that $U\supset\supset V$. The restriction morphism $F(U) \to F(V)$ factors through $\lind {W \supset\supset V}F(W)$. Taking the projective limit we obtain the result. \qed

\begin{oss} The notion of locally weakly quasi-compact can be extended to the case of a site, by generalizing the hypothesis LWC1-LWC3. For our purpose we are interested in the topological setting and we refer to \cite{Phd} for this approach.
\end{oss}

\subsection{c-soft sheaves on locally weakly quasi-compact spaces}


Let $X$ be a locally weakly quasi-compact space, and consider the
category $\mod(k_X)$.

\begin{df} We say that a sheaf $F$ on $X$ is c-soft if the restriction
morphism $\Gamma(W;F) \to \lind{U \supset\supset V}\Gamma(U;F)$ is
surjective for each $V,W \in \op^c(X)$ with $V \subset\subset W$.
\end{df}

It follows from the definition that injective sheaves and flabby
sheaves are c-soft. Moreover, it follows from Corollary
\ref{faclim} that filtrant inductive limits of c-soft sheaves are
c-soft.

\begin{prop}\label{cs1} Let $\exs{F'}{F}{F''}$ be an exact sequence in
$\mod(k_X)$, and assume that $F'$ is c-soft. Then the sequence
$$\exs{\lind{U \supset\supset V}\Gamma(U;F')}{\lind{U \supset\supset V}\Gamma(U;F)}{\lind{U \supset\supset V}\Gamma(U;F'')}$$
is exact for any $V \in \op^c(X)$.
\end{prop}
\dim\ \ Let $s'' \in \lind{U \supset\supset V}\Gamma(U;F'')$. Then
there exists $U \supset\supset V$ such that $s''$ is represented
by $s''_U \in \Gamma(U;F'')$. Let $\{U_i\}_{i \in I} \in \cov(U)$
such that there exists $s_i \in \Gamma(U_i;F)$ whose image is
$s''_U|_{U_i}$ for each $i$. There exists $W \in \op^c(U)$ with
$W\supset\supset V$, a finite covering $\{W_j\}_{j=1}^n$ of $W$
and a map $\varepsilon:J \to I$ of the index sets such that $W_j
\subset\subset U_{\varepsilon(j)}$. We may argue by induction on
$n$. If $n=2$, set $U_i=U_{\varepsilon(i)}$, $i=1,2$. Then
 $(s_1-s_2)|_{U_1 \cap U_2}$ belongs to $\Gamma(U_1 \cap U_2;F')$,
 and  its restriction defines an element of $\lind {W' \supset\supset W_1 \cap
 W_2}\Gamma(W';F')$, hence it extends to $s' \in \Gamma(U;F')$. By
 replacing $s_1$ with $s_1-s'$ on $W_1$ we may assume that
 $s_1=s_2$ on $W_1 \cap W_2$. Then there exists $s \in \Gamma(W_1
 \cup W_2;F)$ with $s|_{W_i}=s_i$. Thus the induction proceeds.
 \qed

\begin{prop}\label{www} Let $\exs{F'}{F}{F''}$ be an exact sequence in
 $\mod(k_{X})$, and assume $F',F$ c-soft. Then $F''$
 is c-soft.
\end{prop}
\dim\ \ Let $V,W \in \op^c(X)$ with $V\subset\subset W$ and let us
consider the diagram below
$$ \xymatrix{\Gamma(W;F) \ar[d]^\alpha \ar[r] & \Gamma(W;F'') \ar[d]^\gamma \\
\lind {U\supset\supset V}\Gamma(U;F) \ar[r]^\beta & \lind {U\supset\supset V}\Gamma(U;F'').} $$  The morphism
$\alpha$ is surjective since $F$ is c-soft and $\beta$ is
surjective by Proposition \ref{cs1}. Then $\gamma$ is
surjective. \qed \\

\begin{prop} The family of c-soft sheaves is injective respect
to the functor $\lind{U \supset\supset V}\Gamma(U;\bullet)$ for each
$V \in \op^c(X)$.
\end{prop}
\dim\ \ The family of c-soft sheaves contains injective sheaves, hence it is cogenerating. Then the result follows from Propositions \ref{cs1} and \ref{www}.
\qed \\

Assume the following hypothesis
\begin{equation}\label{Xcountable}
\text{$X$ has a countable cover $\{U_n\}_{n\in\N}$ with $U_n\in\op^c(X)$, $\forall n\in\N$.}
\end{equation}

\begin{lem}\label{X has a cov} Assume \eqref{Xcountable}. Then there exists a covering $\{V_n\}_{n \in \N}$ of $X$  such
that $V_n \subset \subset V_{n+1}$ and $V_n\in \op ^c(X)$ for each $n \in \N$.
\end{lem}
\dim\ \ Let $\{U_n\}_{n\in\N}$ be a countable cover of $X$ with $U_n \in\op^c(X)$ for each $n\in\N$. Set $V_1=U_1$. Given $\{V_i\}_{i=1}^n$ with $V_{i+1}\supset\supset V_i$, $i=1,\dots,n-1$, let us construct $V_{n+1} \supset\supset V_n$. Consider $x \notin V_n$. Up to take a permutation of $\N$ we may assume $x \in U_{n+1}$. Since $X$ is locally weakly quasi-compact there exists $V_{n+1} \in \op^c(X)$ such that $V_n \cup U_{n+1} \subset\subset V_{n+1}$. \qed

\begin{prop}\label{cs2} Assume \eqref{Xcountable}. Then the category of c-soft sheaves is injective respect
to the functor $\Gamma(X; \bullet)$.
\end{prop}
\dim\ \ Take an exact sequence $\exs{F'}{F}{F''}$ , and suppose
$F'$ c-soft. By Lemma \ref{X has a cov} there exists a covering $\{V_n\}_{n \in \N}$ of $X$  such
that $V_n \subset \subset V_{n+1}$ (and $V_n\in \op ^c(X)$) for each $n \in \N$. All the
sequences
$$\exs{\lind{U_n
\supset\supset V_n}\Gamma(U_n;F')}{\lind{U_n \supset\supset
V_n}\Gamma(U_n;F)}{\lind{U_n \supset\supset V_n}\Gamma(U_n;F'')}$$
are exact by Proposition \ref{cs1}, and the morphism
$\lind{U_{n+1} \supset\supset V_{n+1}}\Gamma(U_{n+1};F') \to
\lind{U_n \supset\supset V_n}\Gamma(U_n;F')$ is surjective for all
$n$. Then by Proposition 1.12.3 of \cite{KS90} the sequence
$$\exs{\lpro n\lind{U_n
\supset\supset V_n}\Gamma(U_n;F')}{\lpro n\lind{U_n \supset\supset
V_n}\Gamma(U_n;F)}{\lpro n\lind{U_n \supset\supset
V_n}\Gamma(U_n;F'')}$$ is exact. By Lemma \ref{lproc} $\lpro n\lind{U_n
\supset\supset V_n}\Gamma(U_n;G) \simeq \Gamma(X;G)$ for any $G
\in \mod(k_X)$ and the result follows. \qed

\begin{es}{\em Let us consider some particular cases


\begin{itemize}
\item[(i)] When $X$ is Noetherian c-soft sheaves are flabby
sheaves.
\item[(ii)] When $X$ has a basis of quasi-compact open subsets, then $F \in \mod(k_X)$ is c-soft if the
restriction morphism $\Gamma(U;F) \to \Gamma(V;F)$ is
surjective, for any quasi-compact open subsets $U,V$ of $X$ with $U \supseteq V$.
\item[(iii)] When $X$ is a locally compact space countable at
infinity, then we recover c-soft sheaves as in chapter II of \cite{KS90}.
\item[(iv)] When $E_{\RP}$ is a vector bundle endowed with the
conic topology, then $F \in \mod(k_{E_{\RP}})$ is c-soft if the
restriction morphism $\Gamma(E_{\RP};F) \to \Gamma(K;F)$ is
surjective, where $K$ is defined as in Example \ref{exiUVi}.
\end{itemize}
}
\end{es}

\section{Sheaves on $\T$-spaces.}\label{Sec tspaces}

In the following we shall assume that $k$ is a field and $X$ is a topological space.
Below we give the definition of $\T$-space, adapting the construction of Kashiwara and Schapira \cite{KS01}. We study the category of sheaves on $X_\T$ generalizing results already known in the case of subanalytic sheaves. Then we prove that as in \cite{ejp} the category of sheaves on $X_\T$ is equivalent to the category of sheaves on a locally weakly-compact topological space $\widetilde{X}_\T$, the $\T$-spectrum, which generalizes the notion of o-minimal spectrum.

\subsection{$\T$-sheaves}

Let $X$ be a topological space and let us consider a family $\T$ of open subsets of $X$.

\begin{df} The topological space $X$ is a $\T$-space if the family $\T$ satisfies the hypotheses below
\begin{equation}\label{hytau}
  \begin{cases}
    \text{(i) $\T$ is a basis for the topology of $X$, and $\varnothing \in \T$}, \\
    \text{(ii) $\T$ is closed under finite unions and intersections},\\
    \text{(iii)  every $U \in \T$ has finitely many $\T$-connected components,}\\

  \end{cases}
\end{equation}
where we define:
\begin{itemize}
\item
 a $\T$-subset is a finite Boolean combination of elements of $\T$;

 \item
 a closed (resp. open) $\T$-subset is a $\T$-subset which is closed (resp. open) in $X$;

 \item
a $\T$-connected subset is a $\T$-subset which is not the disjoint union of two proper $\T$-subsets which are closed and open.
\end{itemize}
\end{df}

\begin{es}\label{T1}
{\em
Let $R=(R,<, 0,1,+,\cdot )$ be a real closed field.
Let $X$ be a locally semialgebraic space (\cite{De91,DK85}) and consider the subfamily of $\op(X)$ defined by $\T= \{U \in \op(X): U \,\, \text{is  semialgebraic}\}$. The family $\T$ satisfy \eqref{hytau}. Note also that the $\T$-subsets of $X$ are exactly the semialgebraic subsets of $X$ (\cite{BCR}).
}
\end{es}

\begin{es}\label{T2}
{\em
Let $X$ be a real analytic manifold and consider the subfamily of $\op(X)$ defined by
$\T=\op^c(X_{sa})=\{U \in \op(X_{sa}): U\,\, $is  subanalytic relatively compact$\}$. The family $\T$ satisfies \eqref{hytau}.
}
\end{es}

\begin{es}\label{T3}
{\em Let $X$ be a real analytic manifold endowed with a subanalytic action $\mu$ of $\RP$. In other words we have a subanalytic map
$$\mu: X \times \RP \to X,$$
which satisfies, for each $t_1,t_2 \in \RP$:
$$
  \begin{cases}
    \mu(x,t_1t_2)=\mu(\mu(x,t_1),t_2), \\
    \mu(x,1)=x.
  \end{cases}
$$
Denote by $X_{\RP}$ the topological space $X$ endowed with the
conic topology,
i.e. $U \in \op(X_{\RP})$ if it is open for the topology of $X$ and invariant by the action of $\RP$. We  will denote by $\op^c(X_{\RP})$ the subcategory of
$\op(X_{\RP})$ consisting of relatively weakly quasi-compact open
subsets.
Consider the subfamily of $\op(X_{\RP})$ defined by
$\T=\op^c(X_{sa,\RP})=\{U \in \op^c(X_{\RP}): U\,\, \text{is  subanalytic}\}$. The family $\T$ satisfies \eqref{hytau}.
}
\end{es}

\begin{es}\label{T4}
{\em
Let ${\mathcal M}=(M,<, (c)_{\in {\mathcal C}}, (f)_{f\in {\mathcal F}}, (R)_{R\in {\mathcal R}} )$ be an arbitrary o-minimal structure.
Let $X$ be a locally definable space (\cite{BO10}) and consider the subfamily of $\op(X)$ defined by  $\T=\op(X_{\rm def})=\{U \in \op(X): U\,\, \text{is  definable}\}$. The family $\T$ satisfies \eqref{hytau}. Note also that (i) the $\T$-subsets of $X$ are exactly the definable subsets of $X$ (by the cell decomposition theorem  in \cite{VD98}, see \cite{ejp} Proposition 2.1).
}
\end{es}

Let $X$ be a $\T$-space. One can endow the category $\T$ with a Gro\-then\-dieck topology, called the $\T$-topology, in the following way: a
family $\{U_i\}_i$ in $\T$ is a covering of $U \in \T$ if it admits a finite subcover. We denote by $X_{\T}$ the associated site,  write for short $k_{\T}$ instead of $k_{X_{\T}}$, and let  $\rho: X \to X_{\T}$ be the natural morphism of sites.
We have functors
\begin{equation}\label{rho}
\xymatrix{\mod(k_X)
\ar@ <2pt> [r]^{\mspace{0mu}\rho_*} &
  \mod(k_\T) \ar@ <2pt> [l]^{\mspace{0mu}\imin \rho}. }
\end{equation}

\begin{prop}\label{fufu} We have $\imin \rho \circ \rho_* \simeq \id$. Equivalently, the functor $\rho_*$ is fully faithful.
\end{prop}
\dim\ \ Let $V \in \op(X)$ and let $G \in \mod(k_\T)$. Then $\imin \rho G = (\rho^\gets F)^{++}$, where $\rho^\gets G \in \psh(k_X)$ is defined by
$$
\op(X) \ni V \mapsto \lind {U \supseteq V, U \in \T}G(U).
$$
In particular, when $U \in \T$, $\rho^\gets G(U)=G(U)$.

Let $F \in \mod(k_X)$ and denote by $\iota:\mod(k_X)\to\psh(k_X)$ the forgetful functor.
The adjunction morphism $\rho^\gets\circ\rho_*\to\id$ in $\psh(k_X)$ defines $\rho^\gets\rho_*F \to \iota F$. This morphism is an isomorphism on $\T$, since $\rho^\gets\rho_*F(U)\simeq\rho_*F(U)\simeq F(U) \simeq \iota F(U)$ when $U \in \T$. By \eqref{hytau} (i) $\T$ forms a basis for the topology of $X$, hence we get an isomorphism
$$
\imin\rho\rho_*F \simeq (\rho^\gets\rho_*F)^{++} \simeq (\iota F)^{++} \simeq F
$$
and the result follows. \qed

\begin{prop}\label{UlimU} Let $\{F_i\}_{i \in I}$ be a filtrant inductive
system in $\mod(k_\T)$ and let $U \in \T$. Then
$$\lind i\Gamma(U;F_i) \iso \Gamma(U;\lind i F_i).$$
\end{prop}
\dim\ \ Denote by $\indl i F_i$ the
presheaf $V \mapsto \lind i \Gamma(V;F_i)$ on $X_\T$. Let $U
\in \T$ and let $S$ be a finite covering of $U$. Since
$\lind i$ commutes with finite projective limits we obtain the
isomorphism $(\indl i F_i)(S) \iso \lind i F_i(S)$ and $F_i(U)
\iso F_i(S)$ since $F_i \in \mod(k_\T)$ for each $i$.
Moreover the family of finite coverings of $U$ is cofinal in
$\cov(U)$. Hence $\indl i F_i \iso (\indl i F_i)^+$. Applying once
again the functor $(\cdot)^+$ we get
$$\indl i F_i \simeq (\indl i F_i)^+ \simeq (\indl i F_i)^{++} \simeq \lind i F_i.$$
Hence applying the functor $\Gamma(U;\cdot)$ we obtain
 the isomorphism $\lind i\Gamma(U;F_i) \iso \Gamma(U;\lind
 i F_i)$ for each $U \in \T.$ \qed \\




\begin{prop}\label{isiuei} Let $F$ be a presheaf on $X_\T$ and
assume that
\begin{itemize}
\item[(i)] $F(\varnothing)=0,$
\item[(ii)] For any $U,V \in \T$ the sequence $\lexs{F(U\cup V)}{F(U) \oplus F(V)}{F(U \cap
V)}$is exact.
\end{itemize}
Then $F \in \mod(k_\T)$.
\end{prop}
\dim\ \ Let $U \in \T$ and let $\{U_j\}_{j=1}^n$ be a
finite covering of $U$. Set for short $U_{ij}=U_i \cap U_j$. We
have to show the exactness of the sequence
$$\lexs{F(U)}{\oplus_{1\leq k\leq n}F(U_k)}{\oplus_{1\leq i<j \leq
n}F(U_{ij})},$$ where the second morphism sends $(s_k)_{1\leq
k\leq n}$ to $(t_{ij})_{1\leq i<j \leq n}$ by
$t_{ij}=s_i|_{U_{ij}}-s_j|_{U_{ij}}$. We shall argue by induction
on $n$. For $n=1$ the result is trivial, and $n=2$ is the
hypothesis. Suppose that the assertion is true for $j \leq n-1$
and set $U'=\bigcup_{1\leq k <n}U_k$. By the induction hypothesis
the following commutative diagram is exact
$$
\xymatrix{&& 0 \ar[d] & 0 \ar[d] \\
0 \ar[r] & F(U) \ar[r] & F(U') \oplus F(U_n) \ar[d] \ar[r]
& F(U'\cap U_n) \ar[d] \\
&& \bigoplus_{i<n}F(U_i) \oplus F(U_n) \ar[d] \ar[r] &
\bigoplus_{i<n}F(U_{in}) \\
&& \bigoplus_{i<j<n}F(U_{ij}).}
$$
Then the result follows. \qed \\

\begin{es} \em{Let us see some examples of sites associated to $\T$-topologies:
\begin{itemize}
\item[(i)] When $\T$ is the family of Example \ref{T1} we obtain the semi-algebraic site of \cite{De91,DK85}.
\item[(ii)] When $\T$ is the family of Example \ref{T2} we obtain the subanalytic site $X_{sa}$ of \cite{KS01,Pr1}.
\item[(iii)] When $\T$ is the family of Example \ref{T3} we obtain the conic subanalytic site of \cite{Pr07}.
\item[(iv)] When $\T$ is the family of Example \ref{T4} we obtain the o-minimal site $X_{\rm def}$. It is the one considered in \cite{ejp} when $X$ is a definable space.
\end{itemize}
}
\end{es}

\subsection{$\T$-coherent sheaves}

Let us consider the category $\mod(k_X)$ of sheaves of
$k_X$-modules on $X$, and denote by $\K$ the subcategory whose
objects are the sheaves $F=\oplus_{i \in I} k_{U_i}$ with $I$
finite and $U_i \in \T$ for each $i$. The following definition is extracted from \cite{KS01}.

\begin{df} Let $\T$ be a subfamily of $\op(X)$ satisfying
$(\ref{hytau})$, and let $F \in \mod(k_X)$.
\begin{itemize}
\item[(i)] $F$ is $\T$-finite if there exists an epimorphism
$G \twoheadrightarrow F$ with $G \in \K$.
\item[(ii)] $F$ is $\T$-pseudo-coherent if for any morphism
$\psi:G \to F$ with $G \in \K$, $\ker \psi$ is $\T$-finite.
\item[(iii)] $F$ is $\T$-coherent if it is both $\T$-finite
and $\T$-pseudo-coherent.
\end{itemize}
\end{df}

Remark that (ii) is equivalent to the same condition with ``$G$ is
$\T$-finite" instead of ``$G \in \K$". One denotes by
$\coh(\T)$ the full subcategory of $\mod(k_X)$ consisting of
$\T$-coherent sheaves. It is easy (see \cite{KS}, Exercise 8.23) to prove that $\coh(\T)$ is additive and stable by kernels.

\begin{lem} \label{lem: finite res} Let $F,G \in \mathcal{K}$. Then, given $\varphi:F \to G$, we have $\ker\varphi \in \mathcal{K}$.
\end{lem}
\dim\ \ We have $F=\oplus_{i=1}^lk_{W_i}$, $G=\oplus_{j=1}^mk_{W'_j}$. Composing with the projection $p_j,\, j=1,..., m$ on each factor of $G$, $\ker\varphi$ will be the intersection of the $\ker p_j\circ \varphi$ so that, if each one has the desired form, the same will happen to their intersection.  Therefore it  is sufficient to assume $m=1$, let us say, $G=k_{W}$. A morphism $\varphi:F \to G$ is then defined by a sequence $v=(v_{1},\dots,v_{l})$, where $v_i$ is the image by $\varphi$ of the section of $k_{W_i}$ defined by $1$ on $W_i$, so $v_{i}=0$ if $W_i\not\subset W$. More precisely, if $s=(s_1,...,s_l)$ is a germ  of $F$ in $y$, we have $\varphi(s_1,..., s_l)=\sum_{i=1}^l{v_i}_{y}{s_i}$. So, given $s=(s_1, ..., s_l)\in\ker\varphi$,  if, for  a given $i$, we have ${v_i}_{y} s_i\neq 0$, then $s$ defines a germ of  $H_i=:\oplus_{i' \neq i}k_{W_{i'}\cap W_{i}}$ in $y$.


Accordingly, $\ker\varphi\simeq\oplus_{i=1}^l H_i$.
\qed \\

Therefore, according to the definition of $\coh(\T)$ and to Lemma \ref{lem: finite res}, any  $F\in\coh(\T)$ admits a finite resolution
$$
K^{\bullet}:=0 \to K^1 \to \cdots \to K^n \to F \to 0
$$
consisting of objects belonging to $\mathcal{K}$.



\begin{prop}\label{rhoUsa} Let $U \in \T$ and consider the constant
sheaf $k_{U_{X_\T}} \in \mod(k_\T)$. We have
$k_{U_{X_\T}} \simeq \rho_*k_U$.
\end{prop}
\dim\ \ Let $F$ be the presheaf on $X_\T$ defined by $F(V)=k$ if
$V \subset U$, $F(V)=0$ otherwise. This is a separated presheaf
and $k_{U_{X_\T}}=F^{++}$. Moreover there is an injective arrow
$F(V) \hookrightarrow \rho_*k_U(V)$ for each $V \in \op(X_\T)$.
Hence $F^{++} \hookrightarrow \rho_*k_U$ since the functor
$(\cdot)^{++}$ is exact. Let $\mathcal{S}\subseteq \T$ be the sub-family of $\T$-connected elements. Then $\mathcal{S}$ forms a basis for
the Grothendieck topology of $X_\T$. For each $W \in \mathcal{S}$ we have $F(W) \simeq
\rho_*k_U(W) \simeq k$ if $W \subset U$ and $F(W)=0$ otherwise.
Then $F^{++} \simeq \rho_*k_U$. \qed \\

\begin{prop} \label{rho exact} The restriction of $\rho_*$ to $\coh(\T)$ is exact.
\end{prop}

\dim\ \
Let us consider an epimorphism
${G}\twoheadrightarrow{F}$ in $\coh(\T)$, we have to
prove that $\psi:{\rho_*G}\to{\rho_*F}$ is an epimorphism. Let $U
\in \T$ and let $0 \neq s \in \Gamma(U;\rho_*F) \simeq
\Ho_{k_X}(k_U,F)$ (by adjunction). Set $G'=G \times_F k_U=\ker(G \oplus k_U \rightrightarrows F)$. Then $G' \in
\coh(\T)$ and moreover $G'\twoheadrightarrow k_U$. There
exists a finite $\{U_i\}_{i \in I} \subset \T$ of  $\T$-connected elements  such that
$\oplus_ik_{U_i}\twoheadrightarrow G'$. The composition $k_{U_i}
\to G' \to k_U$ is given by the multiplication by $a_i \in k$. Set
$I_0=\{k_{U_i};\; a_i \neq 0\}$, we may assume $a_i=1$. We get a
diagram
$$
\xymatrix{\oplus_{i \in I_0}k_{U_i} \ar@{->>} [dr] \ar[r] & G' \ar@{->>} [d] \ar[r] & G \ar@{->>} [d] \\
& k_U \ar[r]^s & F.}
$$
The composition $k_{U_i} \to G' \to G$ defines $t_i \in
\Ho_{k_X}(k_{U_i},G)\simeq \Gamma(U_i;\rho_*G)$. Hence for each $s \in
\Gamma(U;\rho_*F)$ there exists a finite covering $\{U_i\}$ of $U$
and $t_i \in \Gamma(U_i;\rho_*G)$ such that $\psi(t_i)=s|_{U_i}$.
This means that $\psi$ is surjective.
\qed \\

\begin{nt} Since the functor $\rho_*$ is fully faithfull and exact on $\coh(\T)$, we will often identify $\coh(\T)$ with its image in $\mod(k_\T)$ and write $F$ instead of $\rho_*F$ for $F \in \coh(\T)$.
\end{nt}

\begin{teo}\label{teo coh stable}
The following hold:
\begin{itemize}
\item[(i)]
The category $\coh(\T)$ is stable by finite sums, kernels, cokernels and extensions in $\mod(k_\T)$.

\item[(ii)]
The category $\coh(\T)$ is stable by  $\bullet \otimes_{k_\T} \bullet $ in $\mod(k_\T)$.
\end{itemize}
\end{teo}
\dim\ \
(i) The result follows from a general result of homological algebra of \cite{KS96}, Appendix A.1. With the notations of \cite{KS96} let ${\bf P}$ be the set of finite families of elements of $\T$, for ${\mathcal U}=\{U_i\}_{i\in I} \in {\bf P}$ set
$$
L({\mathcal U})=\oplus_ik_{U_i},
$$
for ${\mathcal V}=\{V_j\}_{j\in J} \in {\bf P}$ set
$$
\Ho_{{\bf P}}({\mathcal U},{\mathcal V})=\Ho_{k_\T}(L({\mathcal U}),L({\mathcal V}))=\oplus_i\oplus_j\Ho_{k_\T}(k_{U_i},k_{V_j})
$$
and for $F \in \mod(k_\T)$ set
$$
H({\mathcal U},F)=\Ho_{k_\T}(L({\mathcal U}),F)=\oplus_i\Ho_{k_\T}(k_{U_i},F).
$$
By Proposition A.1 of \cite{KS96} in order to prove (i) it is enough to prove the properties (A.1)-(A.4) below:
\begin{itemize}
\item[(A.1)] For any ${\mathcal U}=\{U_i\} \in {\bf P}$ the functor $H({\mathcal U},\bullet)$ is left exact in $\mod(k_\T)$.
\item[(A.2)] For any morphism $g:{\mathcal V} \to {\mathcal W}$ in ${\bf P}$, there exists a morphism $f:{\mathcal U} \to {\mathcal V}$ in ${\bf P}$ such that ${\mathcal U} \stackrel{f}{\to} {\mathcal V} \stackrel{g}{\to} {\mathcal W}$ is exact.
\item[(A.3)] For any epimorphism $f:F \to G$ in $\mod(k_\T)$, ${\mathcal U} \in {\bf P}$ and $\psi \in H({\mathcal U},G)$, there exists ${\mathcal V} \in {\bf P}$ and an epimorphism $g \in \Ho_{\bf P}({\mathcal V},{\mathcal U})$ and $\varphi \in H({\mathcal V},F)$ such that $\psi \circ g = f \circ \varphi$.
\item[(A.4)] For any ${\mathcal U},{\mathcal V} \in {\bf P}$ and $\psi \in H({\mathcal U},L({\mathcal V}))$ there exists ${\mathcal W} \in {\bf P}$ and an epimorphism $f \in \Ho_{\bf P}({\mathcal W},{\mathcal U})$ and a morphism $g \in \Ho_{\bf P}({\mathcal W},{\mathcal U})$ such that $L(g)=\psi \circ f$ in $\Ho_{k_\T}(L({\mathcal W}),L({\mathcal V}))$.
\end{itemize}
It is easy to check that the axioms (A.1)-(A.4) are satisfied.

(ii) Let $F \in \coh(\T)$. Then $F$ has a resolution
$$
\oplus_{j \in J}k_{U_j} \to \oplus_{i\in I}k_{U_i} \to F \to 0
$$
with $I$ and $J$ finite. Let $V \in \T$. The sequence
$$
\oplus_{j \in J}k_{V \cap U_j} \to \oplus_{i \in T}k_{V \cap U_i} \to F_V \to 0
$$
is exact. Then it follows from (i) that $F_V$ is coherent. Let $G \in \coh(\T)$. The sequence
$$
\oplus_{j \in J}G_{U_j} \to \oplus_{i \in I}G_{U_i} \to G \otimes_{k_\T} F \to 0
$$
is exact. The sheaves $G_{U_i}$ and $G_{U_j}$ are coherent for each $i \in I$ and each $j \in J$. Hence it follows by (i) that $G \otimes_{k_\T} F$ is coherent as required.
\qed

\begin{cor}\label{cor coh stable}
The following hold:
\begin{itemize}
\item[(i)]
The category $\coh(\T)$ is stable by finite sums, kernels, cokernels in $\mod(k_X)$.

\item[(ii)]
The category $\coh(\T)$ is stable by  $\bullet \otimes_{k_X} \bullet $ in $\mod(k_X)$.
\end{itemize}
\end{cor}
\dim\ \ (i) The stability under finite sums and kernels is easy, see \cite{KS}, Exercise 8.23. Let $F,G \in \coh(\T)$ and let $\varphi:F \to G$ be a morphism in $\mod(k_X)$. Then $\rho_*(\varphi)$ is a morphism in $\mod(k_\T)$ and 
$\coker(\rho_*\varphi) \in \coh(\T)$ by Theorem \ref{teo coh stable}. We have
$\coker(\rho_*\varphi) \simeq \rho_*\coker\varphi$ since $\rho_*$ is exact on $\coh(\T)$ by Proposition \ref{rho exact}. Composing with $\imin \rho$ and applying Proposition \ref{fufu} we obtain
$\coker\varphi \in \coh(\T)$.

(ii) The proof of the stability by $\bullet\otimes_{k_X}\bullet$ is similar to that of Theorem \ref{teo coh stable}.
\qed

\begin{teo}\label{eqlambda} (i) Let $G \in \coh(\T)$ and let
$\{F_i\}$ be a filtrant inductive system in $\mod(k_\T)$.
Then we have the isomorphism
$$\lind i \Ho_{k_\T}(\rho_*G,F_i) \iso
\Ho_{k_\T}(\rho_*G,\lind i F_i).$$

(ii) Let $F \in \mod(k_\T)$. There exists a small filtrant
inductive system $\{F_i\}_{i \in I}$ in $\coh(\T)$ such
that $F \simeq \lind i \rho_*F_i$.
\end{teo}
\dim\ \ (i) There exists an exact sequence $\rexs{G_1}{G_0}{G}$
with $G_1,G_0$ finite direct sums of constant sheaves $k_U$ with
$U \in \T$. Since $\rho_*$ is exact on
$\coh(\T)$ and commutes with finite sums, by Proposition
\ref{rhoUsa} we are reduced to prove the isomorphism $\lind i
\Gamma(U;F_i) \iso \Gamma(U;\lind i F_i)$. Then the result follows
from Proposition \ref{UlimU}.

(ii) Let $F \in \mod(k_\T)$, and define
\begin{eqnarray*}
I_0 & := & \{(U,s):\; U \in \T,\; s \in \Gamma(U;F)\} \\
G_0 & := & \oplus_{(U,s) \in I_0}\rho_*k_U
\end{eqnarray*}
The morphism $\rho_*k_U \to F$, where the section $1 \in
\Gamma(U;k_U)$ is sent to $s \in \Gamma(U;F)$ defines un
epimorphism $\varphi:G_0 \to F$. Replacing $F$ by $\ker\varphi$ we
construct a sheaf $G_1=\oplus_{(V,t) \in I_1}\rho_*k_V$ and an
epimorphism $G_1\twoheadrightarrow \ker\varphi$. Hence we get an
exact sequence $\rexs{G_1}{G_0}{F}$. For $J_0\subset I_0$ set for
short $G_{J_0}=\oplus_{(U,s)\in J_0}\rho_*k_U$ and define
similarly $G_{J_1}$. Set
$$
J=\{(J_1,J_0);\; J_k\subset I_k,\; J_k \text{ is finite and } {\rm
im} \varphi|_{G_{J_1}} \subset G_{J_0}\}.
$$
The category $J$ is filtrant and $F \simeq \lind {(J_1,J_0)\in J}
\coker(G_{J_1}\to G_{J_0})$. \qed \\

\begin{cor} Let $G \in \coh(\T)$ and let
$\{F_i\}$ be a filtrant inductive system in $\mod(k_\T)$.
Then we have an isomorphism
$$\lind i \ho_{k_\T}(G,F_i) \iso
\ho_{k_\T}(G,\lind i F_i).$$
\end{cor}
\dim\ \ Let $U \in \T$. We have the chain of isomorphisms
\begin{eqnarray*}
\Gamma(U;\lind i \ho_{k_\T}(G,F_i)) & \simeq & \lind i \Gamma(U;\ho_{k_\T}(G,F_i)) \\
& \simeq & \lind i \Ho_{k_\T}(G_U,F_i)) \\
& \simeq &  \Ho_{k_\T}(G_U,\lind iF_i)) \\
& \simeq &  \Gamma(U;\ho_{k_\T}(G,\lind i F_i)),
\end{eqnarray*}
where the first and the third isomorphism follow from Theorem \ref{eqlambda} (i). the fact that $G_U \in \coh(\T)$ follows from Theorem \ref{teo coh stable} (ii). \qed \\


As in \cite{KS01}, we can define the indization of the category $\coh(\T)$. Recall that the category ${\rm Ind}(\coh(\T))$, of ind-$\T$-coherent sheaves is the category whose objects are filtrant inductive limits of functors
$$
\lind i \Ho_{\coh(\T)}(\bullet,F_i) \ \ \ \ \text{($\indl i F_i$ for short)},
$$
where $F_i \in \coh(\T)$, and the morphisms are the natural transformations of such functors. Note that since $\coh(\T)$ is a small category, ${\rm Ind}(\coh(\T))$ is equivalent to the category of $k$-additive left exact contravariant functors from $\coh(\T)$ to $\mod(k).$  See \cite{KS} for a complete exposition on indizations of categories.
We can extend the functor $\rho_*:\coh(\T)\to \mod(k_{\T})$ to $\lambda:{\rm Ind}(\coh(\T)) \to \mod(k_{\T})$ by setting $\lambda(\indl i F_i):=\lind i \rho_* F_i$.

\begin{cor}
The functor  $\lambda:{\rm Ind}(\coh(\T)) \to \mod(k_{\T})$ is an equivalence of categories.
\end{cor}
\dim \ \ Let $F=\indl j F_j, G=\indl i G_i \in \I(\coh(\T))$. By Theorem \ref{eqlambda} (i) and the fact that  the functor $\rho_*$ is fully faithfull  on $\coh(\T)$ we have
\begin{eqnarray*}
\Ho_{k_\T}(\lambda(F),\lambda(G)) & \simeq & \Ho_{k_\T}(\lind j \rho_* F_j,\lind i \rho_*G_i) \\
& \simeq & \lpro j \lind i \Ho_{k_\T}(\rho_* F_j,\rho_*G_i) \\
& \simeq & \lpro j \lind i \Ho_{\coh(\T)}(F_j,G_i) \\
& \simeq & \Ho_{{\rm Ind}(\coh(\T))}(F,G),
\end{eqnarray*}
hence $\lambda$ is fully faithful. By Theorem \ref{eqlambda} (ii) for each $F \in \mod(k_{\T})$ there exists $G=\indl i F_i \in {\rm Ind}(\coh(\T))$ such that $\lambda(G)=\lind i \rho_*F_i \simeq F$, hence $\lambda$ is essentially surjective.
\qed \\

\subsection{$\T$-flabby sheaves}\label{quinj}

\begin{df}
We say that an object $F \in \mod(k_\T)$ is $\T$-flabby if for each $U,V \in \T$ with $V \supseteq U$
the restriction morphism $\Gamma(V;F) \to \Gamma(U;F)$ is
surjective.
\end{df}


\begin{oss} Remark that the category $\mod(k_\T)$ is a Grothendieck category, hence it has enough injectives.  It follows from the definition that injective sheaves are $\T$-flabby. This implies that the family of $\T$-flabby objects is cogenerating in $\mod(k_\T)$.
\end{oss}

\begin{es} {\em Let us see some examples of $\T$-flabby sheaves:
\begin{itemize}
\item[(i)] When $\T$ is the family of Example \ref{T1} we obtain the family of $sa$-flabby objects of \cite{De91}.
\item[(ii)] When $\T$ is the family of Example \ref{T2} we obtain the family of quasi-injective objects of \cite{Pr1}.
\end{itemize}
}
\end{es}


\begin{prop}
The following hold:
\begin{itemize}
\item[(i)]
Let $F_i$ be a filtrant inductive system of $\T$-flabby sheaves. Then $\lind i F_i$ is $\T$-flabby.

\item[(ii)] Products of $\T$-flabby objects  are $\T$-flabby.
\end{itemize}
\end{prop}
\dim\ \
We will only prove (i) since the proof of (ii) is similar since taking products is exact and commutes with taking sections. Let $U \in \T$. Then for each $i$ the restriction morphism $\Gamma(V;F_i) \to \Gamma(U;F_i)$ is surjective. Applying the exact $\lind i$ and using Proposition \ref{UlimU}, the morphism
$$
\Gamma(V;\lind i F_i) \simeq \lind i \Gamma(V;F_i) \to \lind i \Gamma(U;F_i) \simeq \Gamma(U;\lind i F_i)
$$
is surjective.
\qed

\begin{prop} \label{qinjU}
The full additive subcategory of $\mod (k_\T)$ of $\T$-flabby object is $\Gamma (U;\bullet )$-injective for every $U \in \T$, i.e.:
\begin{itemize}
\item[(i)]
For every $F\in \mod (k_\T)$ there exists a $\T$-flabby object $F'\in \mod (k_\T)$ and an exact sequence $0\rightarrow F\rightarrow F'.$
\item[(ii)]
Let  $\exs{F'}{F}{F''}$ be an exact sequence in
$\mod(k_\T)$ and assume that $F'$ is $\T$-flabby. 
Then the sequence
$$\exs{\Gamma(U;F')}{\Gamma(U;F)}{\Gamma(U;F'')}$$
is exact.
\item[(iii)]
Let $F',F,F'' \in \mod(k_\T)$, and consider the
exact sequence
$$
\exs{F'}{F}{F''}.
$$
Suppose that $F'$ is $\T$-flabby. Then $F$ is $\T$-flabby if and only if $F''$ is $\T$-flabby.
\end{itemize}
\end{prop}
\dim\ \
(i) It follows from the definition that injective sheaves are $\T$-flabby. So (i) holds since it is true for injective sheaves. Indeed, as a Grothendieck category, $\mod(k_\T)$ admits enough injectives.

(ii) Let $s'' \in \Gamma(U;F'')$, and let $\{V_i\}_{i=1}^n\in \cov (U)$ be such that there exists $s_i \in
\Gamma(V_i;F)$ whose image is $s''|_{V_i}$. For $n \geq 2$ on $V_1
\cap V_2$ $s_1-s_2$ defines a section of $\Gamma(V_1 \cap V_2;F')$
which extends to $s' \in \Gamma(U;F')$ since $F'$ is $\T$-flabby. Replace $s_1$ with
$s_1-s'$ (identifying $s'$ with it's image in $F$). We may suppose that $s_1=s_2$ on $V_1 \cap V_2$. Then
there exists $t \in \Gamma(V_1 \cup V_2,F)$ such that
$t|_{V_i}=s_i$, $i=1,2$. Thus the induction proceeds.

(iii) Let $U,V\in \T$ with $V \supseteq U$ and let us
consider the diagram below
$$ \xymatrix{0 \ar[r] & \Gamma(V;F') \ar[d]^\alpha \ar[r] & \Gamma(V;F) \ar[d]^\beta \ar[r] & \Gamma(V;F'') \ar[d]^\gamma \ar[r] & 0\\
0 \ar[r] & \Gamma(U;F') \ar[r] & \Gamma(U;F) \ar[r] & \Gamma(U;F'') \ar[r] & 0} $$
where the row are exact by (ii) 
and the morphism
$\alpha$ is surjective since $F'$ is $\T$-flabby. It follows from the five lemma that $\beta$ is surjective if and only if $\gamma$ is surjective.
\qed \\

\begin{teo}\label{teo flabby hom and Hom}
Let $F \in \mod(k_\T)$. Then the following hold:
\begin{itemize}
\item[(i)]
$F$ is $\T$-flabby if and only if the functor $\Ho_{k_\T}(\bullet ,F)$ is exact on $\coh(\T)$.

\item[(ii)]
If $F$ is $\T$-flabby then  the functor $\ho _{k_\T}(\bullet ,F)$ is exact on $\coh(\T)$.
\end{itemize}

\end{teo}

\dim\ \
(i) is a consequence of a general result of homological algebra  (see Theorem 8.7.2 of \cite{KS}). For (ii),  let $F \in \mod(k_\T)$ be $\T$-flabby. There is
an isomorphism of functors
$$\Gamma(U;\ho _{k_\T}(\bullet ,F)) \simeq \Ho_{k_\T}((\bullet )_U,F)$$
for each $U \in \T$. By Theorem \ref{teo coh stable} and (i)  the
functor $\Ho_{k_\T}((\bullet )_U,F)$
is exact  on $\coh(\T)$ and so the functor $\ho _{k_\T}(\bullet ,F)$ is also exact on $\coh(\T)$. \qed


\begin{teo}\label{rcinj}
Let  $G \in \coh(\T)$. Then the following hold:
\begin{itemize}
\item[(i)]
 The family of $\T$-flabby sheaves is injective with respect to
the functor $\Ho_{k_\T}(G,\bullet )$.

\item[(ii)]
The family of $\T$-flabby sheaves is injective with
respect to the functor $\ho_{k_\T}(G,\bullet )$.
\end{itemize}
\end{teo}

\dim\ \
(i) Let $G \in \coh(\T)$.  Let  $\exs{F'}{F}{F''}$ be an exact sequence in
$\mod(k_\T)$ and assume that $F'$ is $\T$-flabby.
We have to show that the sequence
$$\exs{\Ho_{k_\T}(G,F')}{\Ho_{k_\T}(G,F)}{\Ho_{k_\T}(G,F'')}$$
is exact.

There is an epimorphism $\varphi:\oplus_{i \in I}
k_{U_i} \to G$ where $I$ is finite and $U_i \in \T$ for
each $i \in I$. The sequence $\exs{\ker \varphi}{\oplus_{i \in I} k_{U_i}}{G}$ is
exact. We set for short $G_1=\ker\varphi$ and $G_2=\oplus_{i \in
I} k_{U_i}$. We get the following diagram where the first column is
exact  by Theorem \ref{teo flabby hom and Hom} (i)
$$
\xymatrix{ & 0 \ar[d] & 0 \ar[d] & 0 \ar[d] & \\
0 \ar[r] & \Ho_{k_\T}(G,F') \ar[d] \ar[r] &
\Ho_{k_\T}(G,F)
\ar[d] \ar[r] & \Ho_{k_\T}(G,F'') \ar[d] \ar[r] & 0 \\
0 \ar[r] & \Ho_{k_\T}(G_2,F') \ar[d] \ar[r] &
\Ho_{k_\T}(G_2,F)
\ar[d] \ar[r] & \Ho_{k_\T}(G_2,F'') \ar[d] \ar[r] & 0\\
0 \ar[r] & \Ho_{k_\T}(G_1,F') \ar[d] \ar[r] &
\Ho_{k_\T}(G_1,F)
\ar[d] \ar[r] & \Ho_{k_\T}(G_1,F'') \ar[d] \ar[r] & 0\\
& 0  & 0  & 0  & }
$$

The second row is exact by Proposition \ref{qinjU} (ii), hence the top row is exact by the snake lemma.  \\

(ii) Let $G \in \coh(\T)$. It is enough to check that
for each $U \in \T$ and each exact sequence
$\exs{F'}{F}{F''}$ with $F'$ $\T$-flabby, the
sequence
$$\exs{\Gamma(U;\ho_{k_\T}(G,F'))}{\Gamma(U;\ho_{k_\T}(G,F))}{\Gamma(U;\ho_{k_\T}(G,F''))}$$
is exact. We have
$$\Gamma(U,\ho_{k_\T}(G,\bullet )) \simeq \Ho_{k_\T}(G_U,\bullet ),$$ and,  by (i) and the fact that $G_U\in \coh(\T)$ (Theorem \ref{teo coh stable} (ii)),  $\T$-flabby objects are
injective with respect to the functor
$\Ho_{k_\T}(G_U,\bullet )$ for each $G \in \coh(\T)$,
and for each $U \in \T$.
\qed \\



\begin{prop}\label{hoinj}
Let $F \in \mod(k_\T)$. Then $F$ is
$\T$-flabby if and only if $\ho_{k_\T}(G,F)$ is $\T$-flabby for
each $G \in \coh(\T)$.
\end{prop}
\dim\ \
Suppose that  $F$ is $\T$-flabby, and let $G \in
\coh(\T)$.
We have
$$\Ho_{k_\T}(\bullet ,\ho_{k_\T}(G,F)) \simeq \Ho_{k_\T}(\bullet  \otimes_{k_\T} G,F)$$
and $\Ho_{k_\T}(\bullet \otimes_{k_\T} G,F)$ is exact on $\coh(\T)$ by Theorems \ref{teo coh stable} (ii)  and  \ref{teo flabby hom and Hom} (i).

Suppose that $\ho_{k_\T}(G,F)$ is $\T$-flabby for each $G \in
\coh(\T)$. Let $U,V \in \T$ with $V \supseteq U$. For each $W \in \T$ the morphism $\Gamma(V;\Gamma_WF) \to \Gamma(U;\Gamma_WF)$ is surjective. Hence the morphism
\begin{eqnarray*}
\Gamma(V;F) & \simeq & \Gamma(V;\Gamma_VF) \\
& \to & \Gamma(U;\Gamma_VF) \\
& \simeq & \Gamma(U;F)
\end{eqnarray*}
is surjective. \qed

Let us consider the following subcategory of $\mod(k_\T)$:
$$
\mathcal{P}_{X_\T} := \{G \in \mod(k_\T);\; G \text{ is
$\Ho_{k_\T}(\bullet ,F)$-acyclic for each $F\in\F_{X_\T}$}\},
$$
where ${\mathcal F}_{X_\T}$ is the family of $\T$-flabby objects of $\mod(k_\T)$.

This category is generating. In fact if $\{U_j\}_{j\in J} \in \T$, then $\oplus_{j\in J} k_{U_j} \in \mathcal{P}_{X_\T}$ by Theorem \ref{rcinj} (and the fact that
$$\Pi \Ho _{k_\T}(\bullet , \bullet ) \simeq \Ho _{k_\T}(\oplus \bullet , \bullet )$$
 and products are exact). Moreover $\mathcal{P}_{X_\T}$ is stable by $\bullet
\otimes_{k_\T} K$, where $K \in \coh(\T)$. In fact if $G \in \mathcal{P}_{X_\T}$ and $F \in \F_{X_\T}$ we have
$$\Ho_{k_\T}(G \otimes _{k_\T}K,F) \simeq \Ho_{k_\T}(G,\ho_{k_\T}(K,F))$$
and $\ho_{k_\T}(K,F)$ is $\T$-flabby by Proposition \ref{hoinj}. In particular, if $G\in  \mathcal{P}_{X_\T}$ then $G_U\in  \mathcal{P}_{X_\T}$ for every $U\in \op (X_\T).$

\begin{teo}\label{mathcalJP}
The category $(\mathcal{P}^{op}_{X_\T},\F_{X_\T})$
is injective with respect to
the functors $\Ho_{k_\T}(\bullet ,\bullet )$ and $\ho_{k_\T}(\bullet ,\bullet )$.
\end{teo}

\dim\ \
(i) Let $G \in \mathcal{P}_{X_\T}$ and consider an exact sequence
$\exs{F'}{F}{F''}$ with $F'$ $\T$-flabby. We have to
prove that the sequence
$$\exs{\Ho_{k_\T}(G,F')}{\Ho_{k_\T}(G,F)}{\Ho_{k_\T}(G,F'')}$$
is exact. Since the functor $\Ho_{k_\T}(G,\bullet )$ is acyclic
on $\T$-flabby sheaves we obtain the result.

Let $F$ be $\T$-flabby, and let $\exs{G'}{G}{G''}$ be
an exact sequence on $\mathcal{P}_{X_\T}$. Since the objects of
$\mathcal{P}_{X_\T}$ are $\Ho_{k_\T}(\bullet ,F)$-acyclic the sequence
$$\exs{\Ho_{k_\T}(G'',F)}{\Ho_{k_\T}(G,F)}{\Ho_{k_\T}(G',F)}$$
is exact.


(ii) Let $G \in \mathcal{P}_{X_\T}$, and let
$\exs{F'}{F}{F''}$ be an exact sequence with $F'$ $\T$-flabby. We shall show that for each $U \in
\T$ the sequence
$$\exs{\Gamma(U;\ho_{k_\T}(G,F'))}{\Gamma(U;\ho_{k_\T}(G,F))}{\Gamma(U;\ho_{k_\T}(G,F''))}$$
is exact. This is equivalent  to show that for each $U \in
\T$ the sequence
$$\exs{\Ho_{k_\T}(G_U,F')}{\Ho_{k_\T}(G_U,F)}{\Ho_{k_\T}(G_U,F'')}$$
is exact. This follows since $G_U \in \mathcal{P}_{X_\T}$ as we saw above. The proof of
the exactness in $\mathcal{P}^{op}_{X_\T}$ is similar.
\qed

\begin{prop}\label{prop gamma z and flabby}
Let $F \in \mod(k_\T)$. The following assumptions are equivalent
\begin{itemize}
\item[(i)] $F$ is $\T$-flabby,
\item[(ii)] $F$ is $\Ho_{k_\T}(G,\bullet )$-acyclic for each $G \in \coh(\T)$,
\item[(iii)] $R^1\Ho_{k_\T}(k_{V \setminus U},F)=0$ for each $U,V \in \T$.
\end{itemize}
\end{prop}
\dim\ \
(i) $\Rightarrow$ (ii) follows from Theorem \ref{rcinj}, (ii) $\Rightarrow$ (iii) setting $G=k_{V \setminus U}$  with $U,V \in \T$, (iii) $\Rightarrow$ (i) since if $R^1\Ho_{k_\T}(k_{V \setminus U},F)=0$ for each $U,V \in \T$ with $V \supseteq U$, then the restriction $\Gamma(V;F) \to \Gamma(U;F)$ is surjective. 

\qed



Let $X,Y$ be two topological spaces and let $\T \subset \op(X)$, $\T' \subset \op(Y)$ satisfy \eqref{hytau}. Let $f:X \to Y$ be a continuous map. If $\imin f(\T') \subset \T$ then $f$ defines a morphism of sites $f:X_\T \to Y_{\T'}$.

\begin{prop}\label{prop f star and flabby}
Let $f:X_\T \to Y_{\T'}$ be a morphism of sites. $\T$-flabby sheaves are injective with
respect to the functor $f_*$. The functor $f_*$  sends
$\T$-flabby sheaves to $\T'$-flabby sheaves.
\end{prop}
\dim\ \
Let us consider $V \in \T'$. There is an
isomorphism of functors $\Gamma(V;f_*\bullet ) \simeq \Gamma(\imin
f(V);\bullet )$. It follows from  Proposition \ref{qinjU}  that
$\T$-flabby are injective with respect to the functor
$\Gamma(\imin f(V);\bullet )$ for any $V \in \T'$.

 Let $F$ be $\T$-flabby and let $U,V \in \T'$ with $V \supset U$. Then the morphism
$$\Gamma(V;f_*F) = \Gamma(\imin f(V);F) \to \Gamma(\imin f(U);F) = \Gamma(U;f_*F)$$
is surjective.
\qed

\subsection{$\T$-sheaves on locally weakly quasi-compact spaces}

Assume that $X$ is a locally weakly quasi-compact space.

\begin{lem}\label{UcssVTc} For each $U \in \op^c(X)$ there exists $V \in \T$ such that $U \subset\subset V \subset\subset X$.
\end{lem}
\dim\ \ Since $X$ is locally weakly quasi-compact we may find $W \in \op^c(X)$ such that $U\subset\subset W$. By \eqref{hytau} (i) we may find a covering $\{W_i\}_{i\in I}$ of $X$ with $W_i\in\T$ and $W_i\subset\subset X$ for each $i\in I$. Then there exists a finite family $\{W_j\}_{j=1}^\ell$ whose union $V=\bigcup_{j=1}^\ell W_j$ contains $W$. Then $V \in \T$ and $U\subset\subset V\subset\subset X$. \qed

When $X$ is locally weakly quasi-compact we can construct a left adjoint to the functor $\imin \rho$.

\begin{prop}\label{imineta} Let $F \in \mod(k_\T)$, and let $U \in \op(X)$. Then
$$\Gamma(U;\imin \rho F) \simeq \lpro {V \subset\subset U, V \in
\T} \Gamma(V;F)$$
\end{prop}
\dim\ \ By Theorem \ref{eqlambda} we may assume $F=\lind i
\rho_*F_i$, with $F_i \in \coh(\T)$. Then $\imin \rho F
\simeq \lind i \imin \rho \rho_* F_i \simeq \lind i F_i$. We
have the chain of isomorphisms
$$
\begin{array}{ccccc}
\Gamma(U;\imin \rho F)
 & \hspace{-0.2cm}\simeq & \hspace{-1.3cm}\lpro{V\subset\subset U, V \in \T
}\lind{V\subset\subset W} \hspace{-0.2cm}\Gamma(W;\imin \rho F)
 & \hspace{-0.2cm}\simeq & \hspace{-0.7cm}\lpro{V\subset\subset U, V \in \T
}\lind {V \subset\subset W} \hspace{-0.2cm}\Gamma(W; \lind i \imin \rho\rho_*F_i)\\
 & \hspace{-0.2cm}\simeq & \hspace{-0.8cm}\lpro{V\subset\subset U, V \in
\T}\lind{V\subset\subset W,i}\hspace{-0.25cm}\Gamma(W;\imin \rho
\rho_*F_i)
 & \hspace{-0.2cm}\simeq & \hspace{-0.6cm}\lpro{V\subset\subset U, V \in \T}\lind i \hspace{1mm}\Gamma(V;\imin \rho \rho_*F_i)\\
 & \hspace{-0.2cm}\simeq & \hspace{-0.6cm}\lpro{V\subset\subset U, V \in \T}\lind i \hspace{1mm}\Gamma(V;\rho_*F_i)
 & \hspace{-0.2cm}\simeq & \hspace{-1.8cm}\lpro{V\subset\subset U, V \in \T}\Gamma(V;F),
 \end{array}
 $$
where the first and the fourth isomorphisms follow from Lemma \ref{lproc}, the third isomorphism is a consequence of Corollary
\ref{faclim}, and the last isomorphism follows from Proposition \ref{UlimU}. \qed


\begin{prop}\label{eta!} The functor $\imin \rho$ admits a left adjoint,
denoted by $\rho_!$. It satisfies
\begin{itemize}
\item[(i)] for $F \in \mod(k_X)$ and $U \in \T$, $\rho_! F$ is the sheaf associated to the presheaf $U
\mapsto \lind{U\subset\subset V}
\hspace{-0.2cm}\Gamma(V;F)$, \\
\item[(ii)] For $U \in \op(X)$ one has $\rho_!k_U \simeq \lind{V\subset\subset U, V \in \T
}k_V$.
\end{itemize}
\end{prop}
\dim\ \ Let $\widetilde{F} \in \psh(k_{\T})$ be the presheaf $U
\mapsto \lind {U \subset\subset V} \Gamma(V;F)$, and let $G \in
\mod(k_\T)$. We will construct morphisms
$$\xymatrix{\Ho_{\psh(k_\T)}(\widetilde{F},G) \ar@ <2pt>
[r]^{\xi} &
  \Ho_{k_X}(F,\imin \rho G) \ar@ <2pt> [l]^{\vartheta}}.$$

To define $\xi$, let $\varphi:\widetilde{F} \to G$ and $U \in
\op(X)$. Then the morphism $\xi(\varphi)(U):F(U) \to \imin \rho
G(U)$ is defined as follows

$$F(U) \simeq \lpro{V\subset\subset U, V \in \T
}\lind{V\subset\subset W} F(W) \stackrel{\varphi}{\longrightarrow}
\lpro{V\subset\subset U, V \in \T } G(V) \simeq \imin \rho
G(U).$$

On the other hand, let $\psi:F \to \imin \rho G$ and $U \in \T$.
Then the morphism $\vartheta(\psi)(U):\widetilde{F}(U) \to G(U)$
is defined as follows

$$\widetilde{F}(U) \simeq \lind {U \subset\subset V \in \T} F(V)
\stackrel{\psi}{\longrightarrow} \lind {U \subset\subset V \in
\T} \imin \rho G(V) \to G(U).$$

By construction one can check that the morphism $\xi$ and
$\vartheta$ are inverse to each others. Then (i) follows from the
chain of isomorphisms
$$\Ho_{\psh(k_\T)}(\widetilde{F},G) \simeq
\Ho_{k_{\T}}(\widetilde{F}{}^{++},G) \simeq
\Ho_{k_\T}(\widetilde{F}{}^{++},G).$$

To show (ii), consider the following sequence of isomorphisms
\begin{eqnarray*}
\Ho_{k_\T}(\rho_!k_U,F) & \simeq & \Ho_{k_X}(k_U,\imin \rho F)\\
 & \simeq & \lpro {V \subset\subset U, V \in \T}
 \Ho_{k_\T}(k_V,F)\\
& \simeq & \Ho_{k_\T}(\lind {V \subset\subset U, V \in
\T}k_V,F),
\end{eqnarray*}
where the second isomorphism follows from Proposition
\ref{imineta}. \qed

\begin{prop} The functor $\rho_!$ is exact and commutes with
$\Lind$ and $\otimes$.
\end{prop}
\dim\ \ It follows by adjunction that $\rho_!$ is right exact and
commutes with $\Lind$, so let us show that it is also left exact.
With the notations of Proposition \ref{eta!}, let $F \in
\mod(k_X)$, and let $\widetilde{F} \in \psh(k_{\T})$ be the
presheaf $U \mapsto \lind {U \subset\subset V} \Gamma(V;F)$. Then
$\rho_!F \simeq \widetilde{F}{}^{++}$, and the functors $F \mapsto
\widetilde{F}$ and $G \mapsto G^{++}$ are left exact.\\
Let us show that $\rho_!$ commutes with $\otimes$. Let $F,G \in
\mod(k_X)$, the morphism
$$\lind {U \subset\subset V}F(V) \otimes_k \hspace{-0.1cm} \lind {U
\subset\subset V} G(V) \to \lind {U \subset\subset V} (F(V)
\otimes_k G(V))$$ defines a morphism in $\mod(k_\T)$
$$\rho_!F \otimes_{k_\T} \rho_!G \to \rho_!(F \otimes_{k_X} G)$$
by Proposition \ref{eta!} (i). Since $\rho_!$ commutes with
$\Lind$ we may suppose that $F=k_U$ and $G=k_V$ and the result
follows from Proposition \ref{eta!} (ii).
\qed

\begin{prop} The functor $\rho_!$ is fully faithful. In particular
 one has $\imin \rho \circ \rho_! \simeq \id$. Moreover, for $F \in
 \mod(k_X)$ and $G \in \mod(k_\T)$ one has
$$\imin \rho \ho_{k_\T}(\rho_!F,G) \simeq \ho_{k_X}(F,\imin \rho G).$$
\end{prop}
\dim\ \ For $F,G \in \mod(k_X)$ by adjunction we have
$$\Ho_{k_X}(\imin \rho \rho_!F,G) \simeq \Ho_{k_X}(F,\imin \rho \rho_*G)
\simeq \Ho_{k_X}(F,G).$$ This also implies that $\rho_!$ is fully
faithful, in fact
$$\Ho_{k_\T}(\rho_!F,\rho_!G) \simeq \Ho_{k_X}(F,\imin \rho \rho_!G) \simeq
\Ho_{k_X}(F,G).$$ Now let $K,F \in \mod(k_X)$ and $G \in \mod(k_\T)$, we
have
\begin{eqnarray*}
\Ho_{k_X}(K,\imin \rho \ho_{k_\T}(\rho_! F,G)) & \simeq & \Ho_{k_\T}(\rho_! K,\ho_{k_\T}(\rho_! F,G))\\
 & \simeq & \Ho_{k_\T}(\rho_! K \otimes_{k_\T} \rho_! F,G)\\
 & \simeq & \Ho_{k_\T}(\rho_! (K \otimes_{k_X} F),G)\\
 & \simeq & \Ho_{k_X}(K \otimes_{k_X} F,\imin \rho G)\\
 & \simeq & \Ho_{k_X}(K,\ho_{k_X}(F,\imin \rho G)).
\end{eqnarray*}
\qed

Finally let us consider sheaves of rings in $\mod(k_\T)$. If $\A$
is a sheaf of rings in $\mod(k_X)$, then $\rho_*\A$ and $\rho_!\A$
are sheaves of rings in $\mod(k_\T)$.

Let $\A$ be a sheaf of unitary $k$-algebras on $X$, and let
$\widetilde{\A}\in\psh(k_\T)$ be the presheaf defined by the correspondence $\T
\ni  U \mapsto \lind {U \subset\subset V} \Gamma(V;\A)$. Let $F
\in \psh(k_\T)$, and assume that, for $V \subset U$, with $U,V \in
\T$, the following diagram is commutative:

$$\xymatrix{\Gamma(U;\widetilde{\A}) \otimes_k \Gamma(U;F)
\ar[d] \ar[r] & \Gamma(U;F) \ar[d] \\
\Gamma(V;\widetilde{\A}) \otimes_k \Gamma(V;F) \ar[r] &
\Gamma(V;F).}$$ In this case one says that $F$ is a presheaf of
$\widetilde{\A}$-modules on $\T$.

\begin{prop}\label{eta!mod} Let $\A$ be a sheaf of $k$-algebras on $X$, and let
$F$ be a presheaf of $\widetilde{\A}$-modules on $X_\T$. Then
$F^{++} \in \mod(\rho_!\A)$.
\end{prop}
\dim\ \ Let $U \in \T$, and let $r \in
\lind {U\subset\subset V} \Gamma(V;\A)$. Then $r$ defines a morphism
$\lind{U\subset\subset V} \Gamma(V;\A) \otimes_k \Gamma(W;F) \to \Gamma(W;F)$ for each $W \subseteq U$, $W \in \T$, hence an endomorphism of
$(F^{++})|_{U_{X_\T}}\simeq (F|_{U_{X_\T}})^{++}$. This
morphism defines a morphism of presheaves $\widetilde{\A} \to
{\mathcal E}nd(F^{++})$ and $\widetilde{\A}^{++}\simeq\rho_!\A$ by
Proposition \ref{eta!}. Then $F^{++} \in \mod(\rho_!\A)$. \qed

\begin{prop} Assume that $X$ is locally weakly quasi-compact. Let $F \in \mod(k_\T)$ be $\T$-flabby. Then $\imin \rho F$ is c-soft.
\end{prop}
\dim\ \ Recall that if $U \in
\op(X)$ then $\Gamma(U;\imin \rho F) \simeq \lpro {V
\subset\subset U} \Gamma(V;F)$, where $V \in \T$. Let $W \in \op(X)$, $W \subset\subset X$. It follows from 
Lemma \ref{UcssVTc} that every $U' \supset\supset W$, $U' \in \op(X)$ contains $U \in \T$ such that $U \supset\supset W$. Hence
$$
\lind {U'}\Gamma(U';F) \simeq \lind U\Gamma(U;F),
$$
where $U' \supset\supset W$, $U' \in \op(X)$ and $U \in \T$ such that $U \supset\supset W$.
We have
the chain of isomorphisms
\begin{eqnarray*}
\lind U \Gamma(U;\imin \rho F)
& \simeq & \lind U \lpro {V \subset\subset U} \Gamma(V;F) \\
& \simeq & \lind U \Gamma(U;F)
\end{eqnarray*}
where $U \in \T$, $U \supset\supset W$ and $V \in
\T$. The first isomorphism follows from Proposition \ref{imineta} and second one follows since for each $U \supset\supset W$, $U \in \T$, there exists $V \in \T$ such that $U \supset\supset V \supset\supset W$.

Let $V,W \in \op^c(X)$ with $V \subset\subset W$.
Since $F$ is $\T$-flabby and filtrant inductive limits are
exact, the morphism $\lind {W'} \Gamma(W';\imin \rho F) \simeq \lind {W'} \Gamma(W';F) \to
\lind U \Gamma(U;F) \simeq \lind U \Gamma(U;\imin \rho F)$, where $W',U\in\T$, $W' \supset\supset W$, $U \supset\supset V$, is surjective. Hence
$\Gamma(W;\imin \rho F) \to \lind {U \supset\supset V}\Gamma(U;\imin \rho F)$ is surjective. \qed

\subsection{$\T_{loc}$-sheaves}

Let $X$ be a $\T$-space and let
\begin{equation}\label{Tloc}
\T _{loc}=\{U\in \op(X): U \cap W \in \T\,\,\, \textrm{for every} \,\,\, W \in \T\}.
\end{equation}
Clearly, $\varnothing, X \in \T_{loc}$, $\T\subseteq \T_{loc}$ and $\T_{loc}$ is closed under finite intersections.

\begin{df}
We make the following definitions:
\begin{itemize}
\item
a subset $S$ of $X$ is a $\T_{loc}$-subset if and only if  $S\cap V$ is a $ \T$-subset for every $V\in \T$;

\item
a closed (resp. open) $\T_{loc}$-subset is a $\T_{loc}$-subset which is closed (resp. open) in $X$;

\item
a $\T_{loc}$-connected subset is a $\T_{loc}$-subset which is not the disjoint union of two proper clopen $\T_{loc}$-subsets.

\end{itemize}
\end{df}

Observe that if $\{S_i\}_i$ is a family of $\T_{loc}$-subsets such that $\{i:S_i\cap W\neq \emptyset \}$ is finite for every $W\in \T$, then the union and the intersection of the family $\{S_i\}_i$ is a $\T_{loc}$-subset. Also the complement of a $\T_{loc}$-subset is a $\T_{loc}$-subset. Therefore the $\T_{loc}$-subsets  form a Boolean algebra.\\

\begin{es} \em{ Let us see some examples of $\T_{loc}$ subsets:
\begin{itemize}
\item[(i)] Let $\T$ be the family of Example \ref{T1}. Then the $\T_{loc}$ subsets are the locally semi-algebraic subsets of $X$.
\item[(ii)] Let $\T$ be the family of Example \ref{T2}. Then the $\T_{loc}$ subsets are the subanalytic subsets of $X$.
\item[(iii)] Let $\T$ be the family of Example \ref{T3}. Then the $\T_{loc}$ subsets are the conic subanalytic subsets of $X$.
\item[(iv)] Let $\T$ be the family of Example \ref{T4}. Then the $\T_{loc}$ subsets are the locally definable subsets of $X$.
\end{itemize}
}
\end{es}



One can endow $\T_{loc}$ with a Grothendieck topology in the following way: a
family $\{U_i\}_i$ in $\T_{loc}$ is a covering of $U \in \T_{loc}$ if for any
$V\in\T$, there exists a finite subfamily covering $U \cap V$. We denote by $X_{\T_{loc}}$ the associated site,  write for short $k_{\T_{loc}}$ instead of $k_{X_{\T_{loc}}}$, and let
\[
\xymatrix{
& X \ar[ld]^{\rho _{loc}} \ar[rd]^{\rho } &\\
X_{\T_{loc}} \ar[rr] && X_{\T}
 }
 \]
be the natural morphisms of sites.

\begin{oss}\label{TlocT}
The forgetful functor, induced by the natural morphism of sites $X_{\T_{loc}} \rightarrow X_{\T}$,  gives an equivalence of categories
$$\mod(k_{\T_{loc}}) \iso \mod(k_{\T}).$$
The quasi-inverse to the forgetful functor sends $F \in \mod(k_{\T})$ to $F_{loc}\in\mod(k_{\T_{loc}})$ given by  $F_{loc}(U)=\lpro {V\in\T}F(U \cap V)$  for every $U \in \T_{loc}$.

Therefore, we can and will identify $\mod(k_{\T_{loc}})$ with $\mod(k_{\T})$ and apply the previous results for  $\mod(k_{\T})$ to obtain analogues results for $\mod(k_{\T_{loc}}).$
\end{oss}


Recall that $F \in \mod(k_{\T})$ is $\T$-flabby if the restriction $\Gamma(V;F) \to \Gamma(U;F)$ is surjective for any $U,V \in \T$ with $V \supseteq U$. Assume that
\begin{equation}\label{covT'}
\textrm{$X_{\T_{loc}}$ has a countable cover $\{V_n\}_{n \in \N}$ with $V_n \in \T,$ $\forall n\in\N$.}
\end{equation}

\begin{prop}\label{tflc} Let $F \in \mod(k_{\T})$. Then $F$ is $\T$-flabby if and only if the restriction $\Gamma(X;F) \to \Gamma(U;F)$ is surjective for any $U \in \T_{loc}$.
\end{prop}
\dim\ \ Suppose that $F$ is $\T$-flabby. Consider a covering $\{V_n\}_{n \in \N}$ of $X_{\T_{loc}}$  satisfying \eqref{covT'}.
Set $U_n=U \cap V_n$ and $S_n=V_n \setminus U_n$. All the
sequences
$$
\exs{k_{U_n}}{k_{V_n}}{k_{S_n}}
$$
are exact. Since $F$ is $\T$-flabby the sequence
$$
\exs{\Ho_{k_{\T}}(k_{S_n},F)}{\Ho_{k_{\T}}(k_{V_n},F)}{\Ho_{k_{\T}}(k_{U_n},F)}
$$
is exact. Moreover the morphism $\Ho_{k_{\T}}(k_{S_{n+1}},F) \to \Ho_{k_{\T}}(k_{S_{n}},F)$ is surjective for all $n$ since $S_n=S_{n+1} \cap V_n$ is open in $S_{n+1}$. Then by Proposition 1.12.3 of \cite{KS90} the sequence
$$
\exs{\lpro n\Ho_{k_{\T}}(k_{S_n},F)}{\lpro n\Ho_{k_{\T}}(k_{V_n},F)}{\lpro n\Ho_{k_{\T}}(k_{U_n},F)}
$$
is exact. The result follows since $\lpro n\Gamma(U_n;G) \simeq \Gamma(U;G)$ for any $G \in \mod(k_{\T})$ and $U \in \T_{loc}$. The converse is obvious. \qed

\begin{prop}\label{tflcU} The full additive subcategory of $\mod (k_\T)$ of $\T$-flabby object is $\Gamma (U;\bullet )$-injective for every $U \in \T_{loc}.$
\end{prop}
\dim\ \ Take an exact sequence $\exs{F'}{F}{F''}$, and suppose that
$F'$ is $\T$-flabby. Consider a covering $\{V_n\}_{n \in \N}$ of $X_{\T_{loc}}$ satisfying \eqref{covT'}.
Set $U_n=U \cap V_n$. All the
sequences
$$\exs{\Gamma(U_n;F')}{\Gamma(U_n;F)}{\Gamma(U_n;F'')}$$
are exact by Proposition \ref{qinjU}, and the morphism
$\Gamma(U_{n+1};F') \to
\Gamma(U_n;F')$ is surjective for all
$n$. Then by Proposition 1.12.3 of \cite{KS90} the sequence
$$\exs{\lpro n\Gamma(U_n;F')}{\lpro n\Gamma(U_n;F)}{\lpro n\Gamma(U_n;F'')}$$ is exact. Since $\lpro n\Gamma(U_n;G) \simeq \Gamma(U;G)$ for any $G
\in \mod(k_{\T})$ the result follows. \qed

Let $X,Y$ be two topological spaces and let $\T \subset \op(X)$, $\T' \subset \op(Y)$ satisfy \eqref{hytau}. Let $f:X \to Y$ be a continuous map. If $\imin f(\T'_{loc}) \subseteq \T_{loc}$ then $f$ defines a morphism of sites $f:X_{\T_{loc}} \to Y_{\T'_{loc}}$.

\begin{cor} Let $f:X_{\T_{loc}} \to Y_{\T'_{loc}}$ be a morphism of sites.  $\T$-flabby sheaves are injective with respect to the functor $f_*$. The functor $f_*$  sends $\T$-flabby sheaves to $\T'$-flabby sheaves.
\end{cor}
\dim\ \
Let us consider $V \in \T'_{loc}$. There is an isomorphism of functors $\Gamma(V;f_*\bullet ) \simeq \Gamma(\imin f(V);\bullet )$. It follows from  Proposition \ref{tflcU}   that $\T$-flabby are injective with respect to the functor $\Gamma(\imin f(V);\bullet )$ for any $V \in \T'_{loc}$.

 Let $F$ be $\T$-flabby and let $U,V \in \T'$ with $V \supset U$. Then the morphism
$$\Gamma(V;f_*F) = \Gamma(\imin f(V);F) \to \Gamma(\imin f(U);F) = \Gamma(U;f_*F)$$
is surjective by Proposition \ref{tflc}.
\qed

\begin{oss}\label{Tc} An interesting case is when $X$ is a locally weakly quasi-compact space and there exists $\mathcal{S} \subseteq \op(X)$ with $\T=\{U \in \mathcal{S}: \,U \subset\subset X\}$ satisfying \eqref{hytau}.

Assume that $X$ satisfies \eqref{Xcountable}.
Then $X$ has a covering $\{V_n\}_{n \in \N}$ of $X$  such that $V_n \in \T$ and
$V_n \subset \subset V_{n+1}$ for each $n \in \N$.
By Lemma \ref{X has a cov} we may find a covering $\{U_n\}_{n \in \N}$ of $X$  such that $U_n \in \op^c(X)$ and
$U_n \subset \subset U_{n+1}$ for each $n \in \N$. By Lemma \ref{UcssVTc} for each $n \in \N$ there exists $V_n \in \T$ such that $U_n \subset\subset V_n \subset\subset U_{n+1}$.

In this situation Proposition \ref{tflc} and \ref{tflcU} are satisfied.
\end{oss}




\subsection{$\T$-spectrum}

Let $X$ be a topological space and let $\P(X)$ be the power set of $X$. Consider a subalgebra $\F$ of the power set Boolean algebra $\langle \P(X),\subseteq \rangle$. Then $\F$ is closed under finite unions, intersections and complements. We refer to \cite{Jo82} for an introduction to this subject.

The Boolean algebra $\F$ has an associated topological space, that we denote by $S(\F)$, called its Stone space. The points in $S(\F)$ are the ultrafilters $\alpha$ on $\F$.
The topology on $S(\F)$ is generated by a basis of open and closed sets consisting of all sets of the form
$$
\widetilde{A} = \{\alpha \in S(\F): A \in \alpha \},
$$
where $A\in \F$. The space
$S(\F)$ is a compact totally disconnected Hausdorff space. Moreover, for each $A \in \F$, the subspace $\widetilde{A}$ is Hausdorff and compact. 


\begin{df} 
Let $X$ be a $\T$-space and let $\F$ be the Boolean algebra of $\T_{loc}$-subsets of $X$ (i.e. Boolean combinations of elements of $\T_{loc}$).
The topological space $\widetilde{X}_\T$ is the data of:
\begin{itemize}
\item the points of $S(\F)$ such that $U \in \alpha$ for some $U \in \T$,
\item a basis for the topology is given by  the family of subsets $\{\widetilde{U}:U\in\T\}$.
\end{itemize}
We call $\widetilde{X}_\T$ the $\T$-spectrum of $X$.
\end{df}





With this topology, for $U \in \T$, the set $\widetilde{U}$ is quasi-compact in $\widetilde{X}_\T$ since it is quasi-compact in $S(\F)$. Hence $\widetilde{X}_\T$ is locally weakly quasi-compact with a basis of quasi-compact open subsets given by $\{\widetilde{U}:U\in\T\}$. Note that if $X\in \T$, then $\widetilde{X}_\T=\widetilde{X}$ which is a spectral topological space.

\begin{oss} We may also define $\widetilde{X}_\T$ by means of prime filters of elements of $\T$. This is because $\T$-subsets can be written as finite unions and intersections of $\T$-open and $\T$-closed subsets. In this situation an ultrafilter is determined by the prime filter contained in it.
\end{oss}

\begin{prop}\label{isosites} Let $X$ be a $\T$-space. Then there is an equivalence of categories $\mod(k_\T) \simeq \mod(k_{\widetilde{X}_\T})$.
\end{prop}
\dim\ \ Let us consider the functor
\begin{eqnarray*}
\zeta^t:\T & \to & \op(\widetilde{X}_\T) \\
U & \mapsto & \widetilde{U}.
\end{eqnarray*}
This defines a morphism of sites $\zeta:\widetilde{X}_\T \to X_\T$.  Indeed, if  $V\in \T,$ $S\in \cov (V)$, then $\widetilde{S}=\{\widetilde{V}_i: V_i\in S\} \in \cov(\widetilde{V})$. Let $F \in \mod(k_\T)$ and consider the presheaf $\zeta^\gets F \in \psh(k_{\widetilde{X}_\T})$ defined by $\zeta^\gets F(U) = \lind {U \subseteq \widetilde{V}} F(V)$. In particular, if $U=\widetilde{V}$, $V \in \T$, $\zeta^\gets F(U) \simeq F(V)$. In this case, by Corollary \ref{facKK} we have the isomorphisms
$$
\imin \zeta F(\widetilde{V}) = (\zeta^\gets F)^{++}(\widetilde{V}) \simeq \zeta^\gets F(\widetilde{V}) \simeq F(V).
$$
Then for $V \in \T$ we have
$$
\zeta_*\imin\zeta F(V) \simeq \imin \zeta F(\widetilde{V}) \simeq F(V).
$$
This implies $\zeta_*\circ\imin\zeta\simeq\id$. On the other hand, given $\alpha \in \widetilde{X}_\T$ and $G \in \mod(k_{\widetilde{X}_\T})$,
\begin{eqnarray*}
(\imin \zeta \zeta_* G)_\alpha & \simeq & \lind {\widetilde{U} \ni \alpha, U \in \T} \imin \zeta \zeta_* G(\widetilde{U})\\
& \simeq & \lind {\widetilde{U} \ni \alpha, U \in \T} \zeta_*G(U)\\
& \simeq & \lind {\widetilde{U} \ni \alpha, U \in \T}G(\widetilde{U}) \\
& \simeq  & G_\alpha
\end{eqnarray*}
since $\{\widetilde{U}: U \in \T\}$ forms a basis for the topology of $\widetilde{X}_\T$. This implies $\imin \zeta \circ \zeta_* \simeq \id$. \qed

\begin{es}  {\em Let us see some examples of $\T$-spectra.
\begin{itemize}
\item[(i)] When $\T$ is the family of Example \ref{T1} the $\T$-spectrum $ \widetilde{X}_{\T}$ of $X$ is the semilagebraic spectrum of $X$ (\cite{De91}). When $X$ is semialgebraic, then $\widetilde{X}_{\T} =\widetilde{X}$, the semialgebraic spectrum of $X$ from \cite{CR82}.
\item[(ii)] When $\T$ is the family of Example \ref{T2} the $\T$-spectrum $ \widetilde{X}_{\T}$ of $X$ is the subanalytic spectrum of $X$. The equivalence $\mod(k_{\widetilde{X}_{sa}}) \simeq \mod(k_{X_{sa}})$ was used in \cite{Pr10} to bound the homological dimension of subanalytic sheaves.
\item[(iii)] When $\T$ is the family of Example \ref{T4} the $\T$-spectrum $ \widetilde{X}_{\T}$ of $X$ is the o-minimal spectrum of $X$. When $X$ is a definable space, then $\widetilde{X}_{\T} =\widetilde{X}$, the o-minimal spectrum of $X$ from \cite{Pi88,ejp}.
\end{itemize}
}
\end{es}

\section{Examples}\label{Sec examples}

In this section we recall our main examples of $\T$-sheaves.
Good references on o-minimality are, for example, the book \cite{VD98} by van den Dries and the notes \cite{Co00} by Coste. For semialgebraic geometry relevant to this paper the reader should consult the work by Delfs \cite{De91}, Delfs and Knebusch \cite{DK85} and the book \cite{BCR} by Bochnak, Coste and Roy. For subanalytic geometry we refer to the work \cite{BM88} by Bierstone and Milmann.

\subsection{The semialgebraic site}

Let $R=(R,<, 0,1,+,\cdot )$ be a real closed field.
Let $X$ be a locally semialgebraic space and consider the subfamily of $\op(X)$ defined by $\T= \{U \in \op(X): U \,\, \text{is  semialgebraic}\}$. The family $\T$ satisfies \eqref{hytau} and the associated site $X_{\T}$ is the semialgebraic site on $X$ of \cite{De91,DK85}. Note also that: (i) the $\T$-subsets of $X$ are exactly the semialgebraic subsets of $X$ (\cite{BCR}); (ii) $\T_{loc}= \{U \in \op(X): U \,\, \text{is  locally semialgebraic}\}$ and (iii) the $\T_{loc}$-subsets of $X$ are exactly the locally semialgebraic subsets of $X$ (\cite{DK85}).\\

One can show (using triangulation of semialgebraic sets, as in \cite{KS90}) that the family $\coh(\T)$ corresponds to the family of sheaves which are locally constant on a locally semi-algebraic stratification of $X$.
For each $F \in \mod(k_{\T})$ there exists a filtrant inductive system $\{F_i\}_{i\in I}$ in $\coh(\T)$ such that $F \simeq \lind i \rho_* F_i$.\\

The subcategory of $\T$-flabby sheaves corresponds to the subcategory of $sa$-flabby sheaves of \cite{De91} and it is injective with respect to $\Gamma(U;\bullet)$, $U \in \op(X_{\T})$ and $\Ho_{k_{\T}}(G,\bullet)$, $G \in \coh(\T)$. Our results on $\T$-flabby sheaves generalize those for  $sa$-flabby sheaves from \cite{De91}.\\

We call in this case the $\T$-spectrum $ \widetilde{X}_{\T}$ of $X$ the semialgebraic spectrum of $X$. The  points of $\widetilde{X}_{\T}$ are the ultrafilters $\alpha$ of locally semialgebraic subsets of $X$ such that $U \in \alpha$ for some $U \in \op(X_{\T})$. This is a locally weakly quasi-compact space with basis of quasi-compact open subsets given by $\{\widetilde{U}:U\in \op(X_{\T})\}$ and there is an equivalence of categories $\mod(k_{\T}) \simeq \mod(k_{\widetilde{X}_{\T}})$. When $X$ is semialgebraic, then $\widetilde{X}_{\T} =\widetilde{X}$, the semialgebraic spectrum of $X$ from \cite{CR82}, and there is an equivalence of categories $\mod(k_{\T}) \simeq \mod(k_{\widetilde{X}})$ (\cite{De91}).


\subsection{The subanalytic site}

Let $X$ be a real analytic manifold and consider the subfamily of $\op(X)$ defined by
$\T=\op^c(X_{sa})=\{U \in \op(X_{sa}): U\,\, $is  subanalytic relatively compact$\}$. The family $\T$ satisfies \eqref{hytau} and the associated site $X_{\T}$ is the subanalytic site $X_{sa}$ of \cite{KS01,Pr1}. In this case the $\T_{loc}$-subsets are the subanalytic subsets of $X$.\\

The family $\coh(\T)$ corresponds to the family $\mod_{\rc}^c(k_X)$ of $\R$-constructible sheaves with compact support, and for each $F \in \mod(k_{X_{sa}})$ there exists a filtrant inductive system $\{F_i\}_{i\in I}$ in $\mod_{\rc}^c(k_X)$ such that $F \simeq \lind i \rho_* F_i$.\\

The subcategory of $\T$-flabby sheaves corresponds to quasi-injective sheaves and it is injective with respect to $\Gamma(U;\bullet)$, $U \in \op(X_{sa})$ and $\Ho_{k_{X_{sa}}}(G,\bullet)$, $G \in \mod_{\rc}(k_X)$. \\

We call in this case the $\T$-spectrum $ \widetilde{X}_{\T}$ of $X$ the subanalytic spectrum of $X$ and denote it by  $\widetilde{X}_{sa}$. The points of $\widetilde{X}_{sa}$ are the ultrafilters of subanalytic subsets of $X$ such that $U \in \alpha$ for some $U \in \op^c(X_{sa})$. Then there is an equivalence of categories $\mod(k_{X_{sa}}) \simeq \mod(k_{\widetilde{X}_{sa}})$.

\

Let $U \in \op(X_{sa})$ and denote by $U_{X_{sa}}$ the site with the topology induced by $X_{sa}$. This corresponds to the site $X_{\T}$, where
$\T=\op^c(X_{sa}) \cap U$. In this situation \eqref{hytau} is satisfied.

\subsection{The conic subanalytic site}

Let $X$ be a real analytic manifold endowed with a subanalytic action $\mu$ of $\RP$. In other words we have a subanalytic map
$$\mu: X \times \RP \to X,$$
which satisfies, for each $t_1,t_2 \in \RP$:
$$
  \begin{cases}
    \mu(x,t_1t_2)=\mu(\mu(x,t_1),t_2), \\
    \mu(x,1)=x.
  \end{cases}
$$
Denote by $X_{\RP}$ the topological space $X$ endowed with the
conic topology,
i.e. $U \in \op(X_{\RP})$ if it is open for the topology of $X$ and invariant by the action of $\RP$. We  will denote by $\op^c(X_{\RP})$ the subcategory of
$\op(X_{\RP})$ consisting of relatively weakly quasi-compact open
subsets.\\


Consider the subfamily of $\op(X_{\RP})$ defined by
$\T=\op^c(X_{sa,\RP})=\{U \in \op^c(X_{\RP}): U\,\, \text{is  subanalytic}\}$. The family $\T$ satisfies \eqref{hytau} and the associated site $X_{\T}$ is the conic subanalytic site $X_{sa,\RP}$. In this case the $\T_{loc}$-subsets are the conic subanalytic subsets.\\

Set $\coh(X_{sa,\RP})=\coh(\T)$. For each $F \in \mod(k_{X_{sa,\RP}})$ there exists a filtrant inductive system $\{F_i\}_{i\in I}$ in $\coh(X_{sa,\RP})$ such that $F \simeq \lind i \rho_* F_i$.\\

The subcategory of $\T$-flabby sheaves is injective with respect to $\Gamma(U;\bullet)$, $U \in \op(X_{sa,\RP})$ and $\Ho_{k_{X_{sa,\RP}}}(G,\bullet)$, $G \in \coh(X_{sa,\RP})$.\\


We call in this case the $\T$-spectrum $ \widetilde{X}_{\T}$ of $X$ the conic subanalytic spectrum of $X$ and denote it by  $\widetilde{X}_{sa,\RP}$. The points of $\widetilde{X}_{sa,\RP}$ are the ultrafilters $\alpha$ of conic subanalytic subsets of $X$ such that $U \in \alpha$ for some $U \in \op^c(X_{sa,\RP})$. Then there is an equivalence of categories $\mod(k_{X_{sa,\RP}}) \simeq \mod(k_{\widetilde{X}_{sa,\RP}})$.

\subsection{The o-minimal site}

Let ${\mathcal M}=(M,<, (c)_{\in {\mathcal C}}, (f)_{f\in {\mathcal F}}, (R)_{R\in {\mathcal R}} )$ be an arbitrary o-minimal structure.
Let $X$ be a locally definable space and consider the subfamily of $\op(X)$ defined by  $\T=\op(X_{\rm def})=\{U \in \op(X): U\,\, \text{is  definable}\}$. The family $\T$ satisfies \eqref{hytau} and the associated site $X_{\T}$ is the o-minimal site $X_{\rm def}$ of \cite{ejp}. Note also that: (i) the $\T$-subsets of $X$ are exactly the definable subsets of $X$ (by the cell decomposition theorem  in \cite{VD98}, see \cite{ejp} Proposition 2.1); (ii) $\T_{loc}= \{U \in \op(X): U \,\, \text{is  locally definable}\}$ and (iii) the $\T_{loc}$-subsets of $X$ are exactly the locally definable subsets of $X$.\\

Set $\coh(X_{\rm def})=\coh(\T)$. For each $F \in \mod(k_{X_{\rm def}})$ there exists a filtrant inductive system $\{F_i\}_{i\in I}$ in $\coh(X_{\rm def})$ such that $F \simeq \lind i \rho_* F_i$.\\

The subcategory of $\T$-flabby sheaves (or definably flabby sheaves) is injective with respect to $\Gamma(U;\bullet)$, $U \in \op(X_{\rm def})$ and $\Ho_{k_{X_{\rm def}}}(G,\bullet)$, $G \in \coh(X_{\rm def})$.\\


We call in this case the $\T$-spectrum $ \widetilde{X}_{\T}$ of $X$ the definable or o-minimal spectrum of $X$ and denote it by  $\widetilde{X}_{\rm def}$. The  points of $\widetilde{X}_{\rm def}$ are the ultrafilters $\alpha$ of the Boolean algebra of locally definable subsets of $X$ such that $U \in \alpha$ for some $U \in \op(X_{\rm def})$. This is a locally weakly quasi-compact space with basis of quasi-compact open subsets given by $\{\widetilde{U}:U\in \op(X_{\rm def})\}$ and there is an equivalence of categories $\mod(k_{X_{\rm def}}) \simeq \mod(k_{\widetilde{X}_{\rm def}})$. When $X$ is definable, then $\widetilde{X}_{\rm def} =\widetilde{X}$, the o-minimal spectrum of $X$ from \cite{Pi88,ejp}, and there is an equivalence of categories $\mod(k_{X_{\rm def}}) \simeq \mod(k_{\widetilde{X}})$ (\cite{ejp}).\\


Finally observe that  since locally semialgebraic spaces are locally definable spaces in a real closed field and real closed fields are o-minimal structures and, relatively compact subanalytic sets are definable sets in the o-minimal expansion of the field of real numbers  by restricted globally analytic functions, both the semialgebraic and subanalytic sheaf theory are special cases of the o-minimal sheaf theory.



\end{document}